\documentclass[12pt,reqno]{amsart}
\usepackage{amssymb,comment}
\usepackage{hyperref}
\usepackage{stmaryrd}
\usepackage{slashed}
\usepackage{url}
\usepackage{amsthm}

\newcounter{count}[section]
\renewcommand{\thecount}{\arabic{section}.\arabic{count}}

\numberwithin{equation}{section}

\newtheorem*{lem1}{Lemma A.1}

\newenvironment{Umgeb1rekursiv}[1]
   {\vspace{0.5cm}
     \noindent
     \par\noindent
     \refstepcounter{count}
     \textbf{#1~\thecount}
     \hspace{0.2cm}
     \itshape
   }
   {\vspace{0.2cm}\par}

\newenvironment{Umgeb1}[1]
   {\vspace{0.5cm}
     \par\noindent
     \refstepcounter{count}
     \textbf{#1~\thecount}
     \hspace{0.2cm}
   }
   {\vspace{0.2cm}\par}

\newenvironment{UmgebBeweis}[1]
   {\vspace{0.5cm}
     \par\noindent
     \textbf{Proof.}
     \hspace{0.2cm}
   }
   {\hfill $\square$
   \vspace{0.2cm}}

\newenvironment{sat}{\begin{Umgeb1rekursiv}{Proposition}}
  {\end{Umgeb1rekursiv}}
\newenvironment{lem}{\begin{Umgeb1rekursiv}{Lemma}}
  {\end{Umgeb1rekursiv}}
  
\newenvironment{theorem}{\begin{Umgeb1rekursiv}{Theorem}}
  {\end{Umgeb1rekursiv}}

\newenvironment{bew}{\begin{UmgebBeweis}{Proof}}
   {\end{UmgebBeweis}}

\newenvironment{bem}{\begin{Umgeb1}{Remark}}
   {\end{Umgeb1}}
   
\newcommand{\R}{\mathbb{R}}

\newcommand{\C}{\mathbb{C}}
\newcommand{\N}{\mathbb{N}}
\newcommand{\Z}{\mathbb{Z}}

\newcommand{\id}{\operatorname{Id}}

\renewcommand{\d}{d}

\newcommand{\tr}{\operatorname{tr}}
\newcommand{\Res}{\operatorname{Res}}

\renewcommand{\part}[2]{\frac{\partial #1}{\partial #2}}

\def\st{\stackrel{\text{def}}{=}}

\title[Boundary value problem at conformally Einstein infinity]{The boundary value problem for Laplacian on differential forms 
and conformally Einstein infinity}
\begin{document}

\author[Matthias Fischmann and Petr Somberg]
{Matthias Fischmann and Petr Somberg}

\address{E. \v{C}ech Institute, Mathematical Institute of Charles University,
Sokolovsk\'a 83, Praha 8 - Karl\'{\i}n, Czech Republic}

\thanks{Research supported by grant GA~CR~P201/12/G028.}

\email{fischmann@karlin.mff.cuni.cz, somberg@karlin.mff.cuni.cz}





\begin{abstract}
  We completely resolve the boundary value problem for differential forms 
  and conformally Einstein infinity in terms of the dual Hahn polynomials. Consequently, we 
	produce explicit formulas for the Branson-Gover operators on Einstein manifolds and prove 
	their representation as a product of second order operators. This leads to an explicit 
	description of $Q$-curvature and gauge companion operators on differential forms. 
\end{abstract}
\maketitle

\tableofcontents
\allowdisplaybreaks

\section{Introduction}
   Boundary value problems have ever played an important role in mathematics 
   and physics. A preferred class of boundary value problems is given by a system 
	 of partial differential equations on manifolds with boundary (or, a submanifold)
	 equipped with a geometrical structure.
	 The representative examples are the Laplace and Dirac operators on Riemannian
	 manifolds with boundaries. 
   An intimately related concept is the Poisson transform and boundary (or, submanifold) 
	 asymptotic of a solution for a system of PDEs, cf. \cite{KKMOOT} 
	 for the case related to compactifications of symmetric spaces.

   In \cite{FG3}, Fefferman and Graham initiated a program allowing to regard a 
	 conformal manifold as the conformal infinity of associated Poincar\'e-Einstein 
	 metric. The boundary value problems on the Poincar\'e-Einstein manifolds 
   for eigenvalue type of differential equations (e.g., Laplace, Dirac, etc.) with 
	 prescribed boundary data are referred to as the boundary value 
	 problems for conformal infinity. It is a remarkable fact that solving 
	 such a boundary value problems leads to an algorithmic (or, recursive) construction of 
   a series of conformally covariant differential operators on functions, 
	 spinors and differential forms \cite{GZ, GMP, AG}. Note that these operators 
	 were originally constructed using 
   the ambient metric of Fefferman and Graham and tractor bundles, cf. 
	 \cite{GJMS, HS, BransonGover}, and were soon recognized to encode 
	 interesting geometrical quantities like 
	 Branson's $Q$-curvature \cite{Branson5}
   or holographic deformations of 
   the Yamabe and Dirac operators \cite{Juhl1, Fischmann}.
The study of these conformally covariant differential 
   operators growed rapidly in the last decades and shed some light 
   on their internal structure, see \cite{Juhl1, GoverLatiniWaldron}. The understanding of such structures is 
   very important for further research because of finding explicit formulas is a difficult task, for example see 
   \cite{LeistnerNurowski, FG3, GoverSilhan, FKS}.


   Let $(M,h)$ be an Einstein manifold. The main result of the present article 
      is the complete and explicit solution of the boundary value problem for 
      the Laplace operator acting on differential forms and the conformally 
      Einstein infinity $(M,h)$. More precisely, we 
      reduce the boundary value problem to a rank two matrix valued system of four step 
      recurrence relations for the coefficients of the asymptotic expansion of form Laplace 
      eigenforms. This combinatorial problem can be resolved in terms of generalized 
      hypergeometric functions, closely related to the dual Hahn polynomials. The results 
      analogous to ours 
      were obtained for scalar and spinor fields in \cite{FG3, FKS}. 
      The key property of reducing the boundary value problem for 
			conformally Einstein infinity is the polynomial 
      character of the $1$-parameter family of metrics given by the Poincar\'e-Einstein metric. 
      As an application, we recover explicit formulas for the Branson-Gover and 
      related $Q$-curvature operators on differential forms on Einstein manifolds and reproduce 
			their factorization 
      as a product of second-order differential operators \cite{GoverSilhan}. 
      This factorizations is not new at all, but we think that our proof is more elementary compared to \cite{GoverSilhan}.

\vspace{5pt}

   Let us briefly indicate the content of our article.
       The Section \ref{Recurrence} is combinatorial in its origin with some
	 implications to hypergeometric function theory. We introduce three series 
	 of polynomials $s^{(-)}_m$, $s^{(+)}_m$ and $s^{(1)}_m$ of degree 
        $m\in\N_0$, depending on spectral parameters. Their origin is motivated by 
        the examples given in Subsection \ref{ConformalEinsteinInfinity}.
        We prove that $s^{(\pm)}_m$ satisfy a three step recurrence relation, cf. Proposition \ref{Identities3} and 
        \ref{Identities4}, while $s^{(1)}_m$ turns out to be a linear combination of $s^{(-)}_k$ for $k=m,m-1,m-2,0$, 
        see Theorem \ref{Identities6}. 
        By a cascade of variable changes, we identify $s^{(\pm)}_m$ and $s^{(1)}_m$ as 
        (a linear combination of) generalized hypergeometric functions, in particular 
        $s^{(\pm)}_m$ are given by the dual Hahn polynomials. 
   
	In Section \ref{BoundaryValueProblem}, we briefly recall the 
   boundary value problem for conformal infinity. First of all, we determine in 
   Theorem \ref{FlatCase} its solution in terms of solution operators 
   when the conformal infinity contains the flat metric. Then we prove in 
	Proposition \ref{EinsteinCase} that the
   polynomials $s^{(\pm)}_m$ and $s^{(1)}_m$ mentioned above are the organizing framework 
   to construct those solution operators for our boundary value problem once the conformal infinity contains an Einstein metric. 
  
	In Section \ref{BransonGoverOperators}, we discuss the emergence of
	 Branson-Gover operators in the framework of solution operators. 
	 Furthermore, we present a new proof of
   Theorem \ref{MainTheorem} and Theorem \ref{BGForExceptionals} found in \cite{GoverSilhan} which state that Branson-Gover operators 
   factorize by second-order differential operators. Finally, we discuss 
   explicit formulas for the gauge companion and the $Q$-curvature operators. 
   
	In Appendices \ref{AppendixA} and \ref{AppendixB} we collect some standard notation, results and properties concerning 
	 generalized hypergeometric functions and Poincar\'e-Einstein metrics.

\section{Some combinatorial identities}\label{Recurrence}
  
  In the present section we discuss a class of special polynomials 
	characterized to satisfy certain recurrence relations. 

    Let $y$ be an abstract variable and define the set of polynomials $R_k(y;\alpha)$, 
    \begin{align}\label{eq:Polynomal1}
      R_k(y;\alpha):=\prod_{l=1}^k \big[y-(\alpha-l)(\alpha-l+1)\big],
    \end{align}
    of degree $k\in\N$ and depending on a parameter $\alpha\in\R$. 
    Conventionally, we set $R_0(y;\alpha):=1$. 
    Furthermore, we introduce the polynomials 
    \begin{align}
      s^{(-)}_m(y)&:=\sum_{k=0}^{m}{m\choose k}
        (\tfrac{\beta}{2}-\lambda-m)_{m-k}(-\tfrac{\beta}{2}-m-1)_{m-k} R_k(y;0)\notag\\
      &=:\sum_{k=0}^{m}C^{(-)}_{k}(m)R_k(y;0),\label{eq:Def-s-}\\
      s^{(+)}_m(y)&:=\sum_{k=0}^m {m\choose k}(\tfrac{\beta}{2}-\lambda-m)_{m-k} 
        (-\tfrac{\beta}{2}-m+1)_{m-k}  R_k(y;0)\notag\\
      &  =:\sum_{k=0}^{m}C^{(+)}_{k}(m)R_k(y;0),\label{eq:Def-s+}
    \end{align}
    for $m\in\N_0$ and two parameters $\beta,\lambda\in\C$. 
    Here we already used the notion of Pochhammer symbol, as reviewed
    in Appendix \ref{AppendixA}. 
    \begin{bem}
      We notice that some versions of the polynomials $s^{(\pm)}_m$ have already appeared in \cite[Chapter $7$]{FG3} and \cite{FKS}. They were used 
      to factorize conformal powers of the Laplace and Dirac operator into second order (respectively, first order) differential operators on Einstein manifolds. 
    \end{bem}

    We shall now recall \cite[Equation $7.19$]{FG3} and prove (since it was not presented) basic recurrence 
    relations satisfied by $s^{(-)}_m$ and $s^{(+)}_m$. 
    \begin{sat}\label{Identities3}
      The collection of polynomials $s^{(-)}_m$, $m\in\N$, satisfies the following recurrence relation 
      \begin{align}
        s^{(-)}_m(y)=\big[y+2m(\lambda+m) &+\tfrac{\beta}{2}(\lambda-\tfrac{\beta}{2}-1) \big]s^{(-)}_{m-1}(y)
        \nonumber \\
         &-(m-1)(\lambda+m)(\lambda-\tfrac{\beta}{2}+m-1)(\tfrac{\beta}{2}+m)s^{(-)}_{m-2}(y), \label{eq:recs-}
      \end{align}
      with $s^{(-)}_0(y)=1$, $s^{(-)}_{-1}(y):=0$. 
    \end{sat}
   
    \begin{bew}
      The identity
      \begin{align*}
        R_{k+1}(y;0)=\big[y-k(k+1)\big] R_{k}(y;0)
      \end{align*}
      for $k\in\N$, leads to
      \begin{align*}
        \big[y+2m&(\lambda+m)+\tfrac{\beta}{2}(\lambda-\tfrac{\beta}{2}-1)\big]s^{(-)}_{m-1}(y)\\
        &=\sum_{k=0}^{m-1}C^{(-)}_k(m-1)\big[y-k(k+1)\big]R_k(y;0)\\
        &\quad+\sum_{k=0}^{m-1}C^{(-)}_k(m-1)
          \big[(2m(\lambda+m)+\tfrac{\beta}{2}(\lambda-\tfrac{\beta}{2}-1)+k(k+1)\big]R_k(y;0)\\
        &=\sum_{k=1}^{m}C^{(-)}_{k-1}(m-1)R_k(y;0)\\
        &\quad+\sum_{k=0}^{m-1}C^{(-)}_k(m-1)
           \big[(2m(\lambda+m)+\tfrac{\beta}{2}(\lambda-\tfrac{\beta}{2}-1)+k(k+1)\big]R_k(y;0).
      \end{align*}
      Therefore, it remains to compare the coefficients by $R_k(y;0)$ on both sides of \eqref{eq:recs-}, 
      which is equivalent to the following set of relations among $C^{(-)}_k(m)$:
      \begin{align*}
        C^{(-)}_{m}(m)&=C^{(-)}_{m-1}(m-1),\\
        C^{(-)}_{m-1}(m)&=C^{(-)}_{m-2}(m-1)\\
        &\quad+C^{(-)}_{m-1}(m-1) \big[ (2m(\lambda+m)+\tfrac{\beta}{2}(\lambda-\tfrac{\beta}{2}-1)+m(m-1)\big],\\
        C^{(-)}_{k}(m)&=C^{(-)}_{k-1}(m-1)\\
        &\quad+C^{(-)}_{k}(m-1) \big[(2m(\lambda+m)+\tfrac{\beta}{2}(\lambda-\tfrac{\beta}{2}-1)+k(k+1)\big]\\
        &\quad-C^{(-)}_k(m-2)\big[(m-1)(\lambda+m)(\lambda-\tfrac{\beta}{2}+m-1)(\tfrac{\beta}{2}+m)\big],
      \end{align*}
      for all $k\in\N_0$ such that $k\leq m-2$, and  $C^{(-)}_{-1}(m):=0$ for all $m\in\N_0$. 
     These relations can be easily verified using the identity
      \begin{align*}
        C^{(-)}_r(m-l)=\frac{{ m-l\choose r} (\tfrac{\beta}{2}-\lambda-m+l)_{m-l-r}(-\tfrac{\beta}{2}-m+l-1)_{m-l-r}}{{m\choose k} (\tfrac{\beta}{2}-\lambda-m)_{m-k}(-\tfrac{\beta}{2}-m-1)_{m-k}} C^{(-)}_k(m)
      \end{align*}
     with $l=1,2$ and $r=0, \ldots , m-1$.  This completes the proof.
    \end{bew}

    \begin{sat}\label{Identities4}
      The collection of polynomials $s^{(+)}_m$, $m\in\N$,  satisfies the recurrence relations 
      \begin{align}
        s^{(+)}_m(y)=\big[y &+2(m-1)(\lambda+m-1)
         +\tfrac{\beta}{2}(\lambda-\tfrac{\beta}{2}+1) \big]s^{(+)}_{m-1}(y) \nonumber \\
         &-(m-1)(\lambda+m-2)(\lambda-\tfrac{\beta}{2}+m-1)(\tfrac{\beta}{2}+m-2)s^{(+)}_{m-2}(y),
      \end{align}
       with  $s^{(+)}_0(y)=1$, $s^{(+)}_{-1}(y):=0$. 
    \end{sat}
    \begin{bew}
      It is completely analogous to the proof of the previous proposition. The claim is equivalent to 
      \begin{align*}
        C^{(+)}_{m}(m)&=C^{(+)}_{m-1}(m-1),\\
        C^{(+)}_{m-1}(m)&=C^{(+)}_{m-2}(m-1)\\
        &\quad+C^{(+)}_{m-1}(m-1) \big[ (2(m-1)(\lambda+m-1)+\tfrac{\beta}{2}(\lambda-\tfrac{\beta}{2}+1)+m(m-1)\big],\\
        C^{(+)}_{k}(m)&=C^{(+)}_{k-1}(m-1)\\
        &\quad+C^{(+)}_{k}(m-1) \big[(2(m-1)(\lambda+m-1)+\tfrac{\beta}{2}(\lambda-\tfrac{\beta}{2}+1)+k(k+1)\big]\\
        &\quad-C^{(+)}_k(m-2)\big[(m-1)(\lambda+m-2)(\lambda-\tfrac{\beta}{2}+m-1)(\tfrac{\beta}{2}+m-2)\big],
      \end{align*}
      for $k\in\N_0$ such that  $k\leq m-2$, and $C^{(+)}_{-1}(m):=0$ for all $m\in\N_0$.
     However, these identities hold due to
      \begin{align*}
        C^{(+)}_r(m-l)=\frac{{ m-l\choose r} (\tfrac{\beta}{2}-\lambda-m+l)_{m-l-r}(-\tfrac{\beta}{2}-m+l+1)_{m-l-r}}{{m\choose k} (\tfrac{\beta}{2}-\lambda-m)_{m-k}(-\tfrac{\beta}{2}-m+1)_{m-k}} C^{(+)}_k(m)
      \end{align*}
     with $l=1,2$ and $r=0, \ldots, m-1$. This completes the proof. 
    \end{bew}

    Furthermore, we introduce another set of polynomials:
    \begin{align}\label{eq:Polynomal2}
      R^{(1)}_k(y):&=\sum_{j=1}^k (\tfrac{\beta}{2}-k+j+1)_{2(k-j)}
         \big[y-\tfrac{\beta}{2}(\tfrac{\beta}{2}+1)\big]R_{j-1}(y;k) 
    \end{align}
    of degree $k\in\N$, and set the convention $R^{(1)}_0(y):=0$. 
    The next Theorem states a non-trivial relation among the polynomials defined in 
    Equations \eqref{eq:Polynomal1} and  \eqref{eq:Polynomal2}.
    \begin{theorem}\label{Identities1}
      The polynomials $R_{m}(\cdot;0)$ and $R^{(1)}_m(\cdot)$ are related by
      \begin{align}
        R^{(1)}_{m}(y) &=R_{m}(y;0)-(\tfrac{\beta}{2}-m+1)_{2m},
      \end{align}
      for all $m\in\N_0$.
    \end{theorem}
    \begin{bew}
      The left hand side is a polynomial 
      in $y$ of degree $m$, and so is the right hand side. 
      Hence, it is sufficient to check that both sides of this polynomial 
      identity have the same value at $m+1$ different points. To that aim, we choose 
      the $m$-tuple $i(i-1)$, $i=1,\ldots, m$, 
      of roots of $R_{m}(y;0)$. We note 
      \begin{multline*}
        R^{(1)}_{m}\big(i(i-1)\big)
         =\sum_{j=1}^m (-1)^j(\tfrac{\beta}{2}-m+j+1)_{2(m-j)}\times\\
         \times(\tfrac{\beta}{2}+i)(\tfrac{\beta}{2}-i+1)\prod_{l=1}^{j-1}\big[(m-l+i)(m-l-i+1)\big],
      \end{multline*}
      and the standard combinatorial identities 
      \begin{align*}
        (\tfrac{\beta}{2}-m+j+1)_{2m-2j}
          &=(-1)^j\tfrac{(\tfrac{\beta}{2}-m+1)_{2m}}{(\tfrac{\beta}{2}-m+1)_j(-\tfrac{\beta}{2}-m)_j},\\
        \prod_{l=1}^{j-1}\big[(m-l+i)(m-l-i+1)\big]&=(-m-i+1)_{j-1}(-m+i)_{j-1}
      \end{align*}
      allow to obtain 
      \begin{align*}
        R^{(1)}_{m}\big(i(i-1)\big)=&(\tfrac{\beta}{2}-m+1)_{2m}(\tfrac{\beta}{2}+i)(\tfrac{\beta}{2}-i+1)
          \sum_{j=1}^m\tfrac{(-m-i+1)_{j-1}(-m+i)_{j-1}}{(\tfrac{\beta}{2}-m+1)_j(-\tfrac{\beta}{2}-m)_j}.
      \end{align*}
      Hence, our claim is equivalent to
      \begin{align*}
        \sum_{j=1}^m\frac{(-m-i+1)_{j-1}(-m+i)_{j-1}}{(\tfrac{\beta}{2}-m+1)_{j-1}(-\tfrac{\beta}{2}-m)_{j-1}}
          =-\frac{(\tfrac{\beta}{2}-m+1)(-\tfrac{\beta}{2}-m)}{(\tfrac{\beta}{2}+i)(\tfrac{\beta}{2}-i+1)}
      \end{align*}
      for all $i=1,\ldots,m$. The identity
      \begin{multline*}
        \frac{1}{ (\tfrac{\beta}{2}+i)(\tfrac{\beta}{2}-i+1)}\Bigg[\frac{(-m-i+1)_{j}(-m+i)_{j}}{(\tfrac{\beta}{2}-m+1)_{j-1}(-\tfrac{\beta}{2}-m)_{j-1}}
          -\frac{(-m-i+1)_{j-1}(-m+i)_{j-1}}{(\tfrac{\beta}{2}-m+1)_{j-2}(-\tfrac{\beta}{2}-m)_{j-2}}\Bigg]\\
        =\frac{(-m-i+1)_{j-1}(-m+i)_{j-1}}{(\tfrac{\beta}{2}-m+1)_{j-1}(-\tfrac{\beta}{2}-m)_{j-1}}
      \end{multline*}
      implies that our sum is a telescoping sum, and the only term which survives the summation process is 
      \begin{align*}
        \sum_{j=1}^m\frac{(-m-i+1)_{j-1}(-m+i)_{j-1}}{(\tfrac{\beta}{2}-m+1)_{j-1}(-\tfrac{\beta}{2}-m)_{j-1}}
          &=-\frac{1}{ (\tfrac{\beta}{2}+i)(\tfrac{\beta}{2}-i+1)}\frac{1}{ (\tfrac{\beta}{2}-m+2)_{-1}(-\tfrac{\beta}{2}-m+1)_{-1}}\\
        &=-\frac{(\tfrac{\beta}{2}-m+1)(-\tfrac{\beta}{2}-m)}{(\tfrac{\beta}{2}+i)(\tfrac{\beta}{2}-i+1)},
      \end{align*}
      because $i=1,\ldots,m$.  
      Finally, for the last evaluation point of our polynomials we take $\tfrac{\beta}{2}(\tfrac{\beta}{2}+1)$: 
      \begin{align*}
        R_m(\tfrac{\beta}{2}(\tfrac{\beta}{2}+1);0)& =\prod_{l=1}^{m}(\tfrac{\beta}{2}(\tfrac{\beta}{2}+1)-l(l-1))
        \\
         & =\prod_{l=1}^{m}\big((\tfrac{\beta}{2}+l)(\tfrac{\beta}{2}-l+1)\big)
          =(\tfrac{\beta}{2}-m+1)_{2m},
      \end{align*}
      which is exactly 
     $R^{(1)}_m(\tfrac{\beta}{2}(\tfrac{\beta}{2}+1))+(\tfrac{\beta}{2}-m+1)_{2m}$,
    since $R^{(1)}_m(\tfrac{\beta}{2}(\tfrac{\beta}{2}+1))=0$.
     This completes the proof.
    \end{bew}

    \begin{theorem}\label{Identities2}
      The set of polynomials $R^{(1)}_m$, $m\in\N_0$, satisfies the following recurrence relations
      \begin{align}
        R^{(1)}_{m+1}(y)=\big[y-m(m+1)\big]R^{(1)}_m(y)
          +(\tfrac{\beta}{2}-m+1)_{2m} \big[y-\tfrac{\beta}{2}(\tfrac{\beta}{2}+1)\big],
      \end{align}
      with  $R^{(1)}_{0}(y)=0$. 
    \end{theorem}
    \begin{bew}
      The proof is straightforward. Starting from  
      \begin{align*}
        R^{(1)}_{m+1}(y)&=(\tfrac{\beta}{2}-m+1)_{2m}\big[y-\tfrac{\beta}{2}(\tfrac{\beta}{2}+1)\big]\\
        &\quad+\sum_{j=2}^{m+1} (\tfrac{\beta}{2}-m+j)_{2m+2-2j}
          \times \big[y-\tfrac{\beta}{2}(\tfrac{\beta}{2}+1) \big] R_{j-1}(y;m+1),
      \end{align*}
      we substitute 
      \begin{align*}
        R_{j-1}(y;m+1)=\big[y-m(m+1)\big]R_{j-2}(y;m)
      \end{align*}
      and  shift the summation index.  This yields 
      \begin{align*}
        R^{(1)}_{m+1}(y)&=(\tfrac{\beta}{2}-m+1)_{2m}\big[y-\tfrac{\beta}{2}(\tfrac{\beta}{2}+1)\big]
         +\big[y-m(m+1)\big]R^{(1)}_m(y),
      \end{align*}
      which completes the proof. 
    \end{bew}

    Finally, we introduce the set of polynomials
    \begin{align}
      s^{(1)}_m(y) & :=  \sum_{k=0}^{m}
        {m\choose k}(\lambda-\tfrac{\beta}{2}+k+1)_{m-k}(\tfrac{\beta}{2}+k+1)_{m-k-1}\times\notag\\
      &\quad\quad\times\big[k(\lambda+\beta+2m)+\tfrac{\beta}{2}(\lambda-\beta-2m)\big]  R^{(1)}_k(y)\notag\\
      & =:  \sum_{k=0}^m D^{(1)}_k(m)R^{(1)}_k(y)\label{eq:Def-s1}
    \end{align}
    of degree $m\in\N_0$, and we remark that $s^{(1)}_0(y)=0$. 
    \begin{theorem}\label{Identities6}
      The set of polynomials $s^{(1)}_m$, $m\in\N_{\geq 2}$, satisfies 
      \begin{align}
         s^{(1)}_m(y)&=(\lambda-\beta+2m)s^{(-)}_{m}(y)\notag\\
         &\quad  -2m(\lambda+2m)(\lambda-\tfrac{\beta}{2}+m)s^{(-)}_{m-1}(y)\notag\\
         &\quad +m(m-1)(\lambda+\beta+2m)(\lambda-\tfrac{\beta}{2}+m-1)_2s^{(-)}_{m-2}(y)\notag\\
         &\quad-(\lambda)_m(\lambda-\beta)(\tfrac{\beta}{2})_m s^{(-)}_0(y).\label{eq:s1Recurr}
      \end{align}
       Notice that $s^{(-)}_0(y)=1$.
    \end{theorem}
    \begin{bew}
     By Theorem \ref{Identities1} and Lemma A.1, we have
      \begin{align*}
         s^{(1)}_m(y)&=\sum_{k=0}^m D^{(1)}_k(m)R_k(y;0)-\sum_{k=0}^m D^{(1)}_k(m)(\tfrac{\beta}{2}-k+1)_{2k}\\
         &=\sum_{k=0}^m D^{(1)}_k(m)R_k(y;0)-(\lambda)_m(\lambda-\beta)(\tfrac{\beta}{2})_m.
      \end{align*}
      Comparison of the coefficients in Equation \eqref{eq:s1Recurr} by $R_k(y;0)$ gives 
      \begin{align*}
        D^{(1)}_m(m)&=(\lambda-\beta+2m)C^{(-)}_m(m),\\
        D^{(1)}_{m-1}(m)&=(\lambda-\beta+2m)C^{(-)}_{m-1}(m)
           -2m(\lambda+2m)(\lambda-\tfrac{\beta}{2}+m)C^{(-)}_{m-1}(m-1),\\
        D^{(1)}_{k}(m)&=(\lambda-\beta+2m)C^{(-)}_{k}(m)
           -2m(\lambda+2m)(\lambda-\tfrac{\beta}{2}+m)C^{(-)}_{k}(m-1),\\
        &\quad +m(m-1)(\lambda+\beta+2m)(\lambda-\tfrac{\beta}{2}+m-1)_2 C^{(-)}_k(m-2),
      \end{align*}
      for all $k\in\N_0$ such that $k\leq m-2$. Checking these identities is straightforward and the 
      proof is complete. 
    \end{bew}

  \subsection{Hypergeometric interpretation of combinatorial identities}
    This subsection has no further application in our paper but we think we should shed some 
    light on our defined polynomials $s^{(\pm)}_m$ and $s^{(1)}_m$, see Equations \eqref{eq:Def-s-}, \eqref{eq:Def-s+} and \eqref{eq:Def-s1}, respectively. 
    It will turn out that we can 
    interpret $s^{(\pm)}_m$ as the
    dual Hahn polynomials, cf. Appendix \ref{AppendixA}. 
    Furthermore, it follows from 
		Theorem \ref{Identities6} that $s^{(1)}_m$ can be realized as 
    a linear combination of the
    dual Hahn polynomials with $y$-independent coefficients. 
		This linear combination can be rewritten 
    as a sum of a hypergeometric polynomials of type $(4,3)$ and $(2,1)$ or, equivalently, as a linear 
    combination with $y$-dependent coefficients of two dual Hahn polynomials. 

    Firstly, we will consider $R_k(\cdot;\alpha)$ in the variable $y(y+1)$, 
    \begin{align}\label{eq:ProductIsPochhammer}
      R_k(y(y+1);\alpha)=(-1)^k (-y-\alpha)_k (y+1-\alpha)_k
    \end{align}
    for all $k\in\N$. 

    Secondly, we observe that by standard Pochhammer identities
    our polynomials $s^{(\pm)}_m$ are given by 
    generalized hypergeometric functions of type $(3,2)$:
    \begin{align*}
      s_m^{(-)}\big(y(y+1)\big)&=(\tfrac{\beta}{2}+2)_{m}(\lambda-\tfrac{\beta}{2}+1)_{m}
        \times {}_3F_2 \left[\begin{matrix}-m\;,\; -y\;,\; 1+y\\ \tfrac{\beta}{2}+2\;,\; \lambda-\tfrac{\beta}{2}+1\end{matrix};1  \right],\\ 
      s_m^{(+)}\big(y(y+1)\big)&=(\tfrac{\beta}{2})_{m}(\lambda-\tfrac{\beta}{2}+1)_{m}
        \times {}_3F_2 \left[\begin{matrix}-m\;,\; -y\;,\; 1+y\\ \tfrac{\beta}{2}\;,\; \lambda-\tfrac{\beta}{2}+1\end{matrix};1  \right].
    \end{align*}
    If we choose $\lambda=\tfrac{\beta}{2}-N$ for some $N\in\N$, 
    we can express them as the dual Hahn polynomials $R_m(\lambda(y);a,b,N)$, see Appendix \ref{AppendixA}, 
    \begin{align*}
       s_m^{(-)}\big(y(y+1)\big)&=(\tfrac{\beta}{2}+2)_{m}(\lambda-\tfrac{\beta}{2}+1)_{m}\times
         R_m\big(y(y+1);\tfrac{\beta}{2}+1,-1-\tfrac{\beta}{2},\tfrac{\beta}{2}-\lambda\big),\\
       s_m^{(+)}\big(y(y+1)\big)&=(\tfrac{\beta}{2})_{m}(\lambda-\tfrac{\beta}{2}+1)_{m}\times
         R_m\big(y(y+1);\tfrac{\beta}{2}-1,1-\tfrac{\beta}{2},\tfrac{\beta}{2}-\lambda\big).
    \end{align*}
    
    \begin{bem}
      One can easily realize that there is no hypergeometric
      series representative for $s^{(\pm)}_m(y)$. For example, the
      subleading coefficients in such an expansion do not factorize
      nicely into linear factors. Moreover, the quotients of
      successive coefficients are not rational functions in the
      summation index (which is,
      in fact, the defining property of a hypergeometric series).
    \end{bem}
    Thirdly, as a consequence of Theorem \ref{Identities6}, the polynomial $s^{(1)}_m$ is a linear 
    combination of generalized hypergeometric functions 
    of type $(3,2)$ and $(2,1)$, 
    \begin{align}\label{eq:HyperDecomp}
     s^{(1)}_m\big(y(y+1)\big)&=(\lambda-\tfrac{\beta}{2}+1)_m\Big[(\lambda-\beta+2m)(\tfrac{\beta}{2}+2)_m
       \times {}_3F_2 \left[\begin{matrix}-m\;,\; -y\;,\; 1+y\\ \tfrac{\beta}{2}+2\;,\; \lambda-\tfrac{\beta}{2}+1\end{matrix};1  \right]\notag\\
     &\quad-2m(\lambda+2m)(\tfrac{\beta}{2}+2)_{m-1}
       \times {}_3F_2 \left[\begin{matrix}-m+1\;,\; -y\;,\; 1+y\\ \tfrac{\beta}{2}+2\;,\; \lambda-\tfrac{\beta}{2}+1\end{matrix};1  \right]\notag\\
     &\quad+m(m-1)(\lambda+\beta+2m)(\tfrac{\beta}{2}+2)_{m-2}
       \times {}_3F_2 \left[\begin{matrix}-m+2\;,\; -y\;,\; 1+y\\ \tfrac{\beta}{2}+2\;,\; \lambda-\tfrac{\beta}{2}+1\end{matrix};1  \right]\notag\\
     &\quad-(\lambda-\beta)(\tfrac{\beta}{2})_m 
       \times  {}_2F_1\left[\begin{matrix}-m\;,\; -\tfrac{\beta}{2}+1\\ \lambda-\tfrac{\beta}{2}+1\end{matrix};1  \right] \Big].
    \end{align}
    Notice that the coefficients of the previous linear combination are $y$-independent. 

    We were informed by Christian Krattenthaler 
		that our polynomial $s^{(1)}_m(y(y+1))$ can be organized by the following two 
    expressions based on various generalized hypergeometric functions:
    \begin{sat}
      The set of polynomials $s^{(1)}_m$, for $m\in\N$, has the following descriptions: 
      \begin{align}
         s^{(1)}_m\big(y(y+1)\big)=(\lambda-\tfrac{\beta}{2}+1)_m(\tfrac{\beta}{2})_m\Big[
           &(\lambda-\beta-2m)\times{}_4F_3 \left[\begin{matrix}  1+\gamma  \;,-m\;,\; -y\;,\; 1+y\\ 
           \gamma\;,\;\tfrac{\beta}{2}+1\;,\; \lambda-\tfrac{\beta}{2}+1\end{matrix};1  \right]\notag \\
         &-(\lambda-\beta)\times{}_2F_1\left[\begin{matrix}-m\;,\; -\tfrac{\beta}{2}+1\\ 
           \lambda-\tfrac{\beta}{2}+1\end{matrix};1  \right] \Big],\label{eq:4F3}
      \end{align}
      where $\gamma:=\tfrac{\beta(\lambda-\beta-2m)}{2(\lambda+\beta+2m)}$. Additionally it holds
      \begin{align}
         s^{(1)}_m\big(y(y+1)\big)&=(\lambda-\tfrac{\beta}{2}+1)_m(\tfrac{\beta}{2})_m(\lambda-\beta-2m)
          \times \Big[ {}_3F_2 \left[\begin{matrix}-m\;,\; -y\;,\; 1+y\\ \tfrac{\beta}{2}+1\;,\; \lambda-\tfrac{\beta}{2}+1\end{matrix};1  \right]\notag\\
          &\quad\quad\quad\quad\quad+\tfrac{2m y(y+1)(\lambda+\beta+2m)}{\beta(\tfrac{\beta}{2}+1)(\lambda-\tfrac{\beta}{2}+1)(\lambda-\beta-2m)}\times{}_3F_2 \left[\begin{matrix}1-m\;,\;1-y\;,\; 2+y\\ \tfrac{\beta}{2}+2\;,\; \lambda-\tfrac{\beta}{2}+2\end{matrix};1  \right]\Big]\notag\\
          &\quad-(\lambda-\tfrac{\beta}{2}+1)_m(\tfrac{\beta}{2})_m(\lambda-\beta)\times 
              {}_2F_1\left[\begin{matrix}-m\;,\; -\tfrac{\beta}{2}+1\\ 
         \lambda-\tfrac{\beta}{2}+1\end{matrix};1  \right], \label{eq:Two3F2} 
      \end{align}
     where the coefficients in the linear combination \eqref{eq:Two3F2} 
      are $y$-depended. 
			\end{sat}
      \begin{bew}
         The proof of Equation \eqref{eq:4F3} is based on the elementary identity 
         \begin{align*}
            (\lambda&-\beta+2m)(\tfrac{\beta}{2}+2)_{m}(-m)_k-2m(\lambda+2m)(\tfrac{\beta}{2}+2)_{m-1}(1-m)_k\\
            &\quad+m(m-1)(\lambda+\beta+2m)(\tfrac{\beta}{2}+2)_{m-2}(2-m)_k\\
            &=\tfrac 14 (\beta+2k+2)(-\beta^2+2\beta k+ \lambda \beta+2 \lambda k-2\beta m+4k m)
              (\tfrac{\beta}{2}+2)_{m-2}(-m)_k\\
            &=(\tfrac{\beta}{2})_m(\lambda-\beta-2m)\tfrac{(\gamma+1)_k(-m)_k}{(\gamma)_k(\tfrac{\beta}{2}+1)_k}.
         \end{align*}
         The standard Pochhammer identity 
				$\tfrac{(\gamma+1)_k}{(\gamma)_k}=1+\tfrac{k}{\gamma}$
         allows to decompose our generalized hypergeometric function 
				${}_4F_3$ into two summands, which lead to Equation \eqref{eq:Two3F2}. 
         The proof is complete.          
      \end{bew}

\section{Boundary value problem for conformal infinity}\label{BoundaryValueProblem}

  We start with a brief reminder on the boundary value problem for conformal infinity 
 and the Laplace operator acting on differential forms, for more details we refer to \cite{AG}. 
 Then we proceed to its complete solution in the case when the conformal infinity 
 contains the flat or an Einstein metric. 

  Let $(M,h)$ be a Riemannian oriented manifold of dimension 
  $n\geq 3$. Note that all statements given below extend to the semi-Riemannian setting by 
  careful checking the number of appearances of minuses induced by the signature. The 
  differential $\d:\Omega^p(M)\to\Omega^{p+1}(M)$ has a formal adjoint given by 
  the codifferential $\delta^h=(-1)^{p+1}(\star^h)^{-1}\circ \d\circ \star^h$ 
  when acting on $p+1$-forms. Here we denoted by $\star^h:\Omega^p(M)\to\Omega^{n-p}(M)$ 
  the Hodge operator on $(M,h)$. The form Laplacian 
  \begin{align}
    \Delta^h:=\d\delta^h+\delta^h\d:\Omega^p(M)\to \Omega^p(M)
  \end{align}
  is formally self-adjoint differential operator of second order.  

  Consider the Poincar\'e-Einstein space $(X,g_+)$ associated to $(M,h)$, see Appendix \ref{AppendixB}.
  A differential $p$-form $\omega$ on $X$ uniquely 
  decomposes (note a different convention compared to \cite{AG}) off the boundary as 
  \begin{align}\label{eq:Splitting}
   \omega=\omega^{(+)}+\tfrac{\d r}{r}\wedge \omega^{(-)},
  \end{align}
  for differential forms $\omega^{(+)}\in\Omega^p(X)$ and 
  $\omega^{(-)}\in\Omega^{p-1}(X)$ characterized by trivial 
  contraction with the normal vector field $\partial_r$.
  It is straightforward to verify that the form Laplacian on $X$ acting 
  in the splitting \eqref{eq:Splitting} is 
  \begin{align}\label{eq:FormLaplacian}
    &\Delta^{g_+}=\begin{pmatrix} -(r\partial_r)^2+(n-2p)r\partial_r
       & 2\d \\ 0 & -(r\partial_r)^2+(n+2-2p)r\partial_r\end{pmatrix}\notag\\
    &\quad+\begin{pmatrix}r^2\Delta^{h_r}-r(\star^{h_r})^{-1}[\partial_r,\star^{h_r}]r\partial_r
       &-r\big[\d,(\star^{h_r})^{-1}[\partial_r,\star^{h_r}]\big]\\
    [r\partial_r,r^2\delta^{h_r}] & r^2\Delta^{h_r}-r\partial_r(\star^{h_r})^{-1}[r\partial_r,\star^{h_r}]\end{pmatrix}=:P+P^\prime.
  \end{align}
  Here $\star^{h_r}$ and $\delta^{h_r}$ denote the Hodge operator and 
  codifferential with respect to 
  the $1$-parameter family of metrics $h_r$ on $M$. In the case when $n$ is odd, the 
  form Laplacian can be expanded as a power series around $r=0$, while in even dimensions $n$ 
  there appear additional $\log(r)$-terms  coming from the Poincar\'e-Einstein metric, 
  cf. \cite[Lemma $2.1$]{AG}.  

  For $\lambda\in\C$, we consider the eigenequation 
  \begin{align}\label{eq:eigenequation}
    \Delta^{g_+}\omega=\lambda(n-2p-\lambda)\omega
  \end{align}
  with $\omega$ a $p$-form on $X$. The boundary value problem for conformal infinity consists of 
  finding an asymptotic solution $\omega\in\Omega^p(X)$ of Equation \eqref{eq:eigenequation} with 
  prescribed boundary value $\varphi\in \Omega^p(M)$. The construction of a solution for this boundary
  value problem 
  is algorithmically described in \cite{AG}. For a manifold with general conformal structure $(M,h)$, 
  performing this algorithm is quite complicated due to the 
  complexity in the construction of the Poincar\'e-Einstein metric. As we shall see in next subsections, 
  there is rather explicit solution when the conformal infinity is 
  metrizable by the flat or an Einstein metric. 

  \subsection{Conformally flat metric}
    Let $(M,h)=(\R^n,\langle\cdot,\cdot\rangle_n)$ be the Euclidean space. 
    Then the associated Poincar\'e-Einstein metric can be realized as the hyperbolic metric 
    \begin{align*}
       g_{+}=x_{n+1}^{-2}(\d x^2_{n+1}+\langle\cdot,\cdot\rangle_n)
    \end{align*}
    on the upper half space $\R^{n+1}_{>0}$. 
    Consider the asymptotic expansion of a $p$-form on $\R^{n+1}$, given by 
    \begin{align}\label{eq:FlatFormalAnsatz}
       \omega=x_{n+1}^\lambda\big[\sum_{j\geq 0}x_{n+1}^j \omega^{(+)}_j
         +\tfrac{\d x_{n+1}}{x_{n+1}}\wedge\sum_{j\geq 0}x_{n+1}^j\omega^{(-)}_j    \big]
    \end{align} 
    for $\omega^{(+)}_j\in\Omega^p(\R^n)$ and $\omega^{(-)}_j\in\Omega^{p-1}(\R^n)$. 
    Formally, one can solve 
    Equation \eqref{eq:eigenequation} for a given initial data $\omega^{(+)}_0=\varphi\in\Omega^p(\R^n)$ 
    in terms of the solution operators 
    \begin{align}\label{eq:FlatSolutionOperators}
       \mathcal{T}^{(+)}_{2j}(\lambda)&:\Omega^p(\R^n)\to\Omega^p(\R^n),\notag\\
       \mathcal{T}^{(-)}_{2j}(\lambda)&:\Omega^{p-1}(\R^n)\to\Omega^{p-1}(\R^n),
    \end{align}
    which are $\langle\cdot,\cdot\rangle_n$-natural differential operators with 
    rational polynomial coefficients in $\lambda$ 
    determining $\omega^{(+)}_{2j}=\mathcal{T}^{(+)}_{2j}(\lambda)\varphi$ and 
    $\omega^{(-)}_{2j}=\mathcal{T}^{(-)}_{2j}(\lambda)\varphi$ uniquely for all $j\in\N_0$. 
		Notice that 
    the solution operators turn out to be well-defined for $\lambda\neq \tfrac n2-p-j$ and 
    $\lambda\neq n-2p$, and that by construction holds $\omega^{(\pm)}_{2j-1}=0$ for $j\in\N$. 
    \begin{bem}
      Due to the absence of curvature, the solution operators 
      $\mathcal{T}^{(\pm)}_{2j}(\lambda)$ are given in terms of $\delta (\d\delta)^{j-1}$, 
      $(\delta\d)^j$ and $(\d\delta)^j$. Here $\delta$ denotes the codifferential with respect to the flat metric $\langle\cdot,\cdot\rangle_n$. 
    \end{bem}  
    \begin{theorem}\label{FlatCase}
     Let $(M,h)$ be the Euclidean space $(\R^n,\langle\cdot,\cdot\rangle_n)$, and
		let $j\in\N$ and $\lambda\not=\frac{n}{2}-p-N$ for all $N\in\{1,\ldots ,j\}$. 
		Then 
      \begin{align}
        \mathcal{T}^{(-)}_{2j}(\lambda)
          &=\tfrac{1}{4^{j-1}(j-1)!(\lambda-\tfrac{n-2p}{2}+1)_{j-1}(\lambda-n+2p)}\delta (\d\delta)^{j-1},\nonumber\\
        \mathcal{T}^{(+)}_{2j}(\lambda)&=\tfrac{1}{4^j j! (\lambda-\tfrac{n-2p}{2}+1)_j(\lambda-n+2p)}
          \big[(\lambda-n+2p)(\delta\d)^j
         +(\lambda-n+2p+2j)(\d\delta)^j  \big]
      \end{align}
      for $\mathcal{T}^{(-)}_{0}(\lambda)=0$, $\mathcal{T}^{(+)}_{0}(\lambda)=\id$. 
      Here $\delta$ denotes the codifferential with respect to the flat metric $\langle\cdot,\cdot\rangle_n$.
    \end{theorem}
    \begin{bew}
      For $(M,h)=(\R^n,\langle\cdot,\cdot\rangle_n)$ and $r=x_{n+1}$, 
      $P^\prime$, see Equation \eqref{eq:FormLaplacian}, reduces to  
      \begin{align*}
        P^\prime=x_{n+1}^2\begin{pmatrix}\Delta& 0 \\
          2\delta &\Delta   \end{pmatrix},
      \end{align*}
      where all operators are considered with respect to $h=\langle\cdot,\cdot\rangle_n$. 
      The ansatz \eqref{eq:FlatFormalAnsatz} solves Equation \eqref{eq:eigenequation} iff the following system is 
      satisfied:
      \begin{align*}
        \big[(2j-2)(2\lambda-n+2p+2j-2)+2(\lambda-n+2p+2j-2)\big]\omega^{(-)}_{2j}
          &=\Delta\omega^{(-)}_{2j-2}+2\delta\omega^{(+)}_{2j-2},\\
        2j(2\lambda-n+2p+2j)\omega^{(+)}_{2j}&=\Delta\omega^{(+)}_{2j-2}+2\d \omega^{(-)}_{2j}
      \end{align*}
      for $j\in\N$ and $\omega^{(+)}_0\in\Omega^p(\R^n)$ arbitrary, while $\omega^{(-)}_0=0$. 
      It is now straightforward to check that the solution operators satisfy the recurrence relations 
      and the proof is complete.
    \end{bew}
    
    \begin{bem}\label{ContinuationFlatCase}
      The Theorem \ref{FlatCase} states the existence of solution operators as long as $\lambda\neq \frac n2-p-N$ for some $N\in\N$. 
      By an abuse of notation we will define $\mathcal{T}_j^{\pm}(\lambda)$ through Theorem \ref{FlatCase} for all $\lambda\in\C$. Thus they become 
      meromorphic in $\lambda$. 
    \end{bem}


  \subsection{Conformally Einstein metrics}\label{ConformalEinsteinInfinity}
  Let $(M,h)$ be an Einstein manifold normalized by 
  $Ric(h)=2\lambda(n-1)h$ for some constant $\lambda\in\R$. This implies that the (normalized)  scalar curvature 
  and the Schouten tensor are given by 
  $J=n\lambda$ and $P=\tfrac{J}{n}h$, 
  respectively.\footnote{We use the convention that $J\st\frac{scal}{2(n-1)}$, where $scal$ denotes the scalar curvature with respect to $h$, and $P\st\frac{1}{n-2}(Ric(h)-J h)$.} 
  In this case the 
  Poincar\'e-Einstein metric is of the form
  \begin{align}\label{eq:EinsteinMetric}
   g_+=r^{-2}\big(\d r^2+J(r)^2 h\big),
  \end{align}
  for $J(r):=(1-\tfrac{J}{2n}r^2)$. 
  From now on we use the abbreviation 
  \begin{align}
    \beta:=n-2p.
  \end{align}
  The explicit formula of the metric \eqref{eq:EinsteinMetric} implies that one can 
  explicitly compute the form Laplacian $\Delta^{g_+}$, especially the term $P^\prime$ in Equation \eqref{eq:FormLaplacian}. 
  \begin{lem}
    Let $(M,h)$ be an Einstein manifold, normalized by $Ric(h)=2\lambda(n-1)h$. Then   
    in the splitting \eqref{eq:Splitting}, it holds  
    \begin{align}
      \Delta^{g_+}=\begin{pmatrix} -(r\partial_r)^2+\beta r\partial_r& 0 \\ 
                   0 & -(r\partial_r)^2+(\beta+2)r\partial_r\end{pmatrix}+
                 \begin{pmatrix} A_1 &A_2\\ A_3 & A_4   \end{pmatrix},
    \end{align}
    where
    \begin{align*}
      A_1=&r^2J(r)^{-2}\Delta^h+\beta \tfrac Jn r^2J(r)^{-1} r\partial_r  ,\\
      A_2=&2(1+\tfrac{J}{2n} r^2)J(r)^{-1}d,\\
      A_3=&2r^2(1+\tfrac{J}{2n}r^2)J(r)^{-3}\delta^h,\\
      A_4=&r^2 J(r)^{-2}\Delta^h+(\beta+2)\tfrac Jn r^2 J(r)^{-2}(2+J(r)r\partial_r).
    \end{align*} 
  \end{lem}
  \begin{bew}
    The explict formula for $h_r$, cf. \eqref{eq:EinsteinMetric}, leads to the explicit form 
    of Equation \eqref{eq:FormLaplacian}. In more detail, we note 
    \begin{align*}
      \star^{h_r}(\eta)&=J(r)^{\beta}\star^h(\eta),
    \end{align*}
    for $\eta\in\Omega^p(M)$. Hence, we get on $p$-forms 
    \begin{align*}
      \delta^{h_r}&=(-1)^{p}(\star_{p-1}^{h_r})^{-1}\circ d\circ \star_p^{h_r}\\
      &=(-1)^{p}J(r)^{2(p-1)-n}J(r)^{n-2p}(\star_{p-1}^h)^{-1}\circ d\circ \star_p^h
        =J(r)^{-2}\delta^h, 
    \end{align*}
    which implies $\Delta^{h_r}=J(r)^{-2}\Delta^h$. 
    Furthermore, for $p$ and $(p-1)$-forms $\omega^{(+)},\omega^{(-)}$ as introduced in 
		\eqref{eq:Splitting}, we have 
    \begin{align*}
      r(\star^{h_r})^{-1}[\partial_r,\star^{h_r}]r\partial_r\omega^{(+)}
        &=-\beta\tfrac Jn J(r)^{-1}r^3\partial_r\omega^{(+)},\\
      r\big[d,(\star^{h_r})^{-1}[\partial_r,\star^{h_r}]\big]\omega^{(-)}&=-2\tfrac Jn r^2J(r)^{-1}d\omega^{(-)},\\
      [r\partial_r,r^2\delta^{h_r}]\omega^{(+)}
        &=2r^2(1+\tfrac{J}{2n}r^2)J(r)^{-3}\delta^h\omega^{(+)},\\
      r\partial_r(\star^{h_r})^{-1}[r\partial_r,\star^{h_r}]\omega^{(-)}
        &=-(\beta+2)r^2\tfrac{J}{n}J(r)^{-2}(2+J(r)r\partial_r)\omega^{(-)}.
    \end{align*}
    The result then follows from \eqref{eq:FormLaplacian}.
  \end{bew}

  The eigenequation \eqref{eq:eigenequation}, acting in the 
  splitting $\omega\simeq \begin{pmatrix} \omega^{(+)}\\ \omega^{(-)}\end{pmatrix}$, 
  is equivalent to the system
  \begin{align}
    \lambda&(\beta-\lambda)\omega^{(-)}+(r\partial_r)^2\omega^{(-)}-(\beta+2)r\partial_r\omega^{(-)}\notag\\
    &=2r^2(1+\tfrac{J}{2n}r^2)J(r)^{-3}\delta^h\omega^{(+)}+r^2 J(r)^{-2}\Delta^h\omega^{(-)}\notag\\
    &\quad+(\beta+2)\tfrac{2J}{n} r^2 J(r)^{-2}\omega^{(-)}
      +(\beta+2)\tfrac{J}{n}r^2 J(r)^{-1}r\partial_r\omega^{(-)}\label{eq:SecondRecurrence},\\
    \lambda&(\beta-\lambda)\omega^{(+)}+(r\partial_r)^2\omega^{(+)}-\beta r\partial_r\omega^{(+)}\notag\\
    &=2(1+\tfrac{J}{2n}r^2)J(r)^{-1}d\omega^{(-)}+r^2J(r)^{-2}\Delta^h\omega^{(+)}
      +\beta\tfrac Jn r^2 J(r)^{-1}r\partial_r\omega^{(+)}\label{eq:FirstRecurrence}.
  \end{align}
  Using the polynomial type of $J(r)=1-\tfrac{J}{2n}r^2$, we multiply 
  Equation \eqref{eq:SecondRecurrence} by $J(r)^3$ and Equation \eqref{eq:FirstRecurrence} 
  by $J(r)^2$. 
  As a result, the coefficients in both Equations \eqref{eq:SecondRecurrence} and \eqref{eq:FirstRecurrence} 
  become polynomials of degree $3$ and $2$ 
  in $r^2$, respectively. This is the key step to formulate an equivalence between a (formal) 
  solution of Equation \eqref{eq:eigenequation} and a solution of a system of recurrence relations given in \eqref{eq:MinusRecur} and \eqref{eq:PlusRecur}. 
  \begin{theorem}\label{EinsteinMainRecurrence}
    Let $(M,h)$ be an Einstein manifold with the normalization 
    given by $Ric(h)=\tfrac{2J}{n}(n-1)h$ for a
    constant $J\in\R$. The eigenequation \eqref{eq:eigenequation} acting on
    \begin{align*}
      \omega=r^{\lambda}\big[\sum_{j\geq 0} r^{j}\omega_j^{(+)}
       +\tfrac{dr}{r}\wedge\sum_{j\geq 0}r^{j}\omega_j^{(-)}\big],
    \end{align*}
    for some (unknown) differential forms 
    $\omega^{(+)}_j\in\Omega^p(M)$ and $\omega^{(-)}_j\in\Omega^{p-1}(M)$, is 
    equivalent to the following recurrence relations:
    \begin{align}
      A_{j} w^{(-)}_{j}&=B_{j}(\tfrac{J}{2n}) w^{(-)}_{j-2}+C_{j} (\tfrac{J}{2n})^2 w^{(-)}_{j-4}
       +D_{j} (\tfrac{J}{2n})^3 w^{(-)}_{j-6}+\Delta^h w^{(-)}_{j-2}-\tfrac{J}{2n}\Delta^h \omega^{(-)}_{j-4}\notag\\
      &\quad+2\delta^h w^{(+)}_{j-2}+2(\tfrac{J}{2n}) \delta^h w^{(+)}_{j-4},\label{eq:MinusRecur}\\
      a_{j} w^{(+)}_{j}&=b_{j}(\tfrac{J}{2n}) w^{(+)}_{j-2}+c_{j} (\tfrac{J}{2n})^2 w^{(+)}_{j-4}
        +\Delta^h w^{(+)}_{j-2}+2d w^{(-)}_{j}-2(\tfrac{J}{2n})^2 dw^{(-)}_{j-4},\label{eq:PlusRecur}
    \end{align}
    with coefficients
    \begin{align}\label{eq:Coefficients}
      A_{j}&:=-\big[(j-2)(\beta-2\lambda-j+2)+2(\beta-\lambda-j+2)\big],\notag\\
      B_{j}&:=-3\big[(j-4)(\beta-2\lambda-j+4)+2(\beta-\lambda-j+4)\big]
        +2(\beta+2)(\lambda+j),\notag\\
      C_{j}&:=3\big[(j-6)(\beta-2\lambda-j+6)+2(\beta-\lambda-j+6)\big]
        -4(\beta+2)(\lambda+j-3),\notag\\
      D_{j}&:=-\big[(j-8)(\beta-2\lambda-j+8)+2(\beta-\lambda-j+8)\big]
        +2(\beta+2)(\lambda+j-6),\notag\\
      a_{j}&:=-j(\beta-2\lambda-j),\notag\\
      b_{j}&:=2\big[(j-2)(\lambda+j-2)+\lambda(\beta+j-2)\big],\notag\\
      c_{j}&:=-(2\lambda+j-4)(\beta+j-4),
    \end{align}
    depending on $j\in\N$, $\beta=n-2p$, $\lambda\in\C$, and  
    the initial data $\omega^{(-)}_0=0$ and $\omega^{(+)}_0\in\Omega^{(p)}(M)$. Furthermore, it holds 
    $\omega^{(\pm)}_{2j-1}=0$ for all $j\in\N$.
  \end{theorem}
  
  \begin{bew}
    Due to the polynomial type of the coefficients in Equations \eqref{eq:SecondRecurrence} and 
    \eqref{eq:FirstRecurrence} after multiplication with appropriate powers of $J(r)$, and the ansatz for 
    $\omega$, we obtain the recurrence relations by 
    comparing the coefficients by $r^{\lambda+j}$ for $j\in\N_0$. Based on  
    $\omega^{(-)}_0=0$ and evenness of involved coefficients in $r$, we 
    get $\omega^{(\pm)}_{2j-1}=0$ for all $j\in\N$. This completes the proof. 
  \end{bew}
	
   In order to get an insight into the solution structure of the recurrence relations \eqref{eq:MinusRecur} 
   and \eqref{eq:PlusRecur}, we present several low-order approximations:
	
   {\bf Second-order approximation:} 
    The relation \eqref{eq:MinusRecur} for $j=2$ gives 
    \begin{align*}
      2(\lambda-\beta)\omega^{(-)}_2=2\delta^h \omega^{(+)}_0,
    \end{align*}
    and \eqref{eq:PlusRecur} for $j=2$ implies 
    \begin{align*}
      4(2\lambda-\beta+2)(\lambda-\beta)\omega^{(+)}_2&
        = 2(\lambda-\beta)R^{(1)}_1\omega^{(+)}_0
        +2(\lambda-\beta+2)R^{(2)}_1\omega^{(+)}_0\\
      &\quad+2(2\lambda-\beta+2)(\lambda-\beta)\beta \tfrac{J}{2n} \omega^{(+)}_0
    \end{align*}
    for
    \begin{align*}
      R_1^{(1)}:=\delta^h d+\tfrac{\beta}{2}(\tfrac{\beta}{2}-1)\tfrac{2J}{n}, \quad R^{(2)}_1:=d\delta^h.
    \end{align*}
    Hence $\omega^{(-)}_2$ and $\omega^{(+)}_2$ are well-defined 
    for $\lambda\neq\beta$ and $\lambda\neq\tfrac{\beta}{2}-1,\beta$, respectively. 
    
  {\bf Fourth-order approximation:} 
    The relation \eqref{eq:MinusRecur} for $j=4$ gives 
    \begin{align*}
      4(2\lambda-\beta+2)(\lambda-\beta)\omega^{(-)}_4
      &=2 R^{(0)}_1\delta^h\omega^{(+)}_0+2(2\lambda-\beta+2)(\beta+4)(\tfrac{J}{2n})\delta^h\omega^{(+)}_0
    \end{align*}
    in terms of the operator
    \begin{align*}
      R^{(0)}_1:=\delta^h d+\tfrac{\beta}{2}(\tfrac{\beta}{2}+1)\tfrac{2J}{n}.
    \end{align*}
 
    The relation \eqref{eq:PlusRecur} for $j=4$ reduces to 
    \begin{align*}
      16(2\lambda-\beta+4)&(2\lambda-\beta+2)(\lambda-\beta)\omega^{(+)}_4\notag\\
      &=2(\lambda-\beta)R^{(1)}_2\omega^{(+)}_0 +2(\lambda-\beta+4)R^{(2)}_2\omega^{(+)}_0\notag\\
      &\quad+4(2\lambda-\beta+4)(\lambda-\beta)(\beta+2)(\tfrac{J}{2n})R^{(1)}_1\omega^{(+)}_0\notag\\
      &\quad+4(2\lambda-\beta+4)\big[2(\lambda-\beta+4)+\beta(\lambda-\beta)\big](\tfrac{J}{2n})R^{(2)}_1\omega^{(+)}_0\notag\\
      &\quad+2(2\lambda-\beta+4)(2\lambda-\beta+2)(\lambda-\beta)\beta(\beta+2)(\tfrac{J}{2n})^2\omega^{(+)}_0,
    \end{align*}
    in terms of additional operators 
    \begin{align*}
      R^{(1)}_2&:=\big[\delta^hd+(\tfrac{\beta}{2}-1)\tfrac{\beta}{2} \tfrac{2J}{n}\big]
       \big[\delta^hd+(\tfrac{\beta}{2}-2)(\tfrac{\beta}{2}+1)\tfrac{2J}{n}\big],\\
      R^{(2)}_2&:= \tfrac{\beta}{2}(\tfrac{\beta}{2}+1)(\tfrac{2J}{n}) d\delta^h
        +d\delta^h\big[d\delta^h+ (\tfrac{\beta}{2}-1)(\tfrac{\beta}{2}+2)\tfrac{2J}{n} \big].
    \end{align*}
    Hence $\omega^{(-)}_4$ and $\omega^{(+)}_4$ are well-defined 
    for $\lambda\neq\tfrac{\beta}{2}-1,\beta$ and $\lambda\neq\tfrac{\beta}{2}-1,\tfrac{\beta}{2}-2,\beta$, 
    respectively. 
  
  {\bf Sixth-order approximation:} 
    The relation \eqref{eq:MinusRecur} for $j=6$ gives 
    \begin{align*}
      16(2\lambda-\beta+2)&(2\lambda-\beta+4)(\lambda-\beta)\omega^{(-)}_6\\
      &=2 R^{(0)}_2 \delta^h\omega^{(+)}_0 
         +4(2\lambda-\beta+4)(\beta+6)(\tfrac{J}{2n}) R^{(0)}_1 \delta^h\omega^{(+)}_0\\
      &\quad+2(2\lambda-\beta+2)(2\lambda-\beta+4)
        (\beta+4)(\beta+6)(\tfrac{J}{2n})^2\delta^h\omega^{(+)}_0
    \end{align*}
    in terms of the additional operator
    \begin{align*}
      R^{(0)}_2:=\big[\delta^h d+\tfrac{\beta}{2}(\tfrac{\beta}{2}+1)\tfrac{2J}{n}\big]
       \big[\delta^h d +(\tfrac{\beta}{2}-1)(\tfrac{\beta}{2}+2)\tfrac{2J}{n}   \big].
    \end{align*}
    The relation \eqref{eq:PlusRecur} for $j=6$ yields after some computations 
    \begin{align*}
      9&6(2\lambda-\beta+6)(2\lambda-\beta+4)(2\lambda-\beta+2)(\lambda-\beta)\omega^{(+)}_6\notag\\
      &=2(\lambda-\beta)R^{(1)}_3\omega^{(+)}_0
        +2(\lambda-\beta+6)R^{(2)}_3\omega^{(+)}_0\notag\\
      &\quad+6(2\lambda-\beta+6)(\lambda-\beta)(\beta+4)(\tfrac{J}{2n}) R^{(1)}_2\omega^{(+)}_0\notag\\
      &\quad+6(2\lambda-\beta+6)\big[4(\lambda-\beta+6)+\beta(\lambda-\beta+2)\big]
          (\tfrac{J}{2n})R^{(2)}_2\omega^{(+)}_0\notag\\
      &\quad+6(2\lambda-\beta+6)(2\lambda-\beta+4)
          (\lambda-\beta)(\beta+2)(\beta+4)(\tfrac{J}{2n})^2 R^{(1)}_1\omega^{(+)}_0\notag\\
      &\quad+6(2\lambda-\beta+6)(2\lambda-\beta+4)\big[2(\lambda-\beta+6)+\beta(\lambda-\beta-2)\big](\beta+4)
          (\tfrac{J}{2n})^2R^{(2)}_1\omega^{(+)}_0\notag\\
      &\quad+2(2\lambda-\beta+6)(2\lambda-\beta+4)(2\lambda-\beta+2)(\lambda-\beta)\beta(\beta+2)(\beta+4)
         (\tfrac{J}{2n})^3\omega^{(+)}_0
    \end{align*}
    expressed in terms of additional operators 
    \begin{align*}
      R^{(1)}_3&:=\big[\delta^hd+(\tfrac{\beta}{2}-1)\tfrac{\beta}{2} \tfrac{2J}{n}\big]
        \big[\delta^hd+(\tfrac{\beta}{2}-2)(\tfrac{\beta}{2}+1)\tfrac{2J}{n}\big]
        \big[\delta^hd+(\tfrac{\beta}{2}-3)(\tfrac{\beta}{2}+2)\tfrac{2J}{n}\big],\\
      R^{(2)}_3&:=(\tfrac{\beta}{2}-1)\tfrac{\beta}{2}(\tfrac{\beta}{2}+1)(\tfrac{\beta}{2}+2)(\tfrac{2J}{n})^2d\delta^h
        +\tfrac{\beta}{2}(\tfrac{\beta}{2}+1)\tfrac{2J}{n}  d\delta^h\big[d\delta^h+ (\tfrac{\beta}{2}-2)(\tfrac{\beta}{2}+3)\tfrac{2J}{n} \big]\\
      &\quad\;+ d\delta^h \big[d\delta^h+(\tfrac{\beta}{2}-1)(\tfrac{\beta}{2}+2)\tfrac{2J}{n}\big]
        \big[d\delta^h + (\tfrac{\beta}{2}-2)(\tfrac{\beta}{2}+3)\tfrac{2J}{n}\big]. 
    \end{align*}
    Hence $\omega^{(-)}_6$ and $\omega^{(+)}_6$ are well-defined for 
    $\lambda\neq\tfrac{\beta}{2}-1,\tfrac{\beta}{2}-2,\beta$ 
    and $\lambda\neq\tfrac{\beta}{2}-1,\tfrac{\beta}{2}-2,\tfrac{\beta}{2}-3,\beta$, respectively. 
    
    \vspace{0.5cm}
   
    The previous approximations indicate the following definition of the solution operators: 
    \begin{align}\label{eq:SolutionOperators}
      \mathcal{T}^{(-)}_m(\lambda)
        :=&\big[(\lambda-\beta) \prod_{k=1}^{m-1}a_{2k}\big]^{-1} (\tfrac{2J}{n})^{m-1}
         s^{(-)}_{m-1}(y^{(-)})\circ \delta^h,\notag\\
      \mathcal{T}^{(+)}_m(\lambda)
        :=&\big[(\lambda-\beta)\prod_{k=1}^{m}a_{2k}\big]^{-1}  (\tfrac{2J}{n})^{m}\bigg[(\lambda-\beta)s^{(+)}_m(y^{(+)})
         +s^{(1)}_m(y^{(1)}) \bigg],
    \end{align}
    for $m\in\N$ and the coefficients $a_k$ introduced in \eqref{eq:Coefficients}. 
		Here we have taken the evaluation of the polynomials $s^{(\pm)}_k$ and $s^{(1)}_k$ at 
    \begin{align*}
      y^{(-)}:=\tfrac{n}{2J}  \delta^h\d+\tfrac{\beta}{2}(\tfrac{\beta}{2}+1),\\
      y^{(+)}:=\tfrac{n}{2J} \delta^h\d+\tfrac{\beta}{2}(\tfrac{\beta}{2}-1),\\
      y^{(1)}:=\tfrac{n}{2J} \d\delta^h+\tfrac{\beta}{2}(\tfrac{\beta}{2}+1).
    \end{align*}
    \begin{bem}
       The inspiration for the definition of 
			$\mathcal{T}^{(+)}_m(\lambda)$ comes from 
			 the scalar case, cf. \cite[Chapter $7$]{FG3}. Restricting to $0$-forms, 
       Equation \eqref{eq:MinusRecur} 
       becomes trivial and $\mathcal{T}^{(-)}_m(\lambda)=0$, while 
			 Equation \eqref{eq:PlusRecur} is solved 
       by $\mathcal{T}^{(+)}_m(\lambda)$ with vanishing term $s^{(1)}_m$.
    \end{bem}
    The proof of the next theorem is mainly based on the combinatorial identities 
    discussed in Section \ref{Recurrence}. 
    \begin{theorem}\label{EinsteinCase}
      Let $m\in\N$, $\varphi\in\Omega^p(M)$ and $\lambda\neq \tfrac{\beta}{2}-N$
      for all $N\in\{1,\ldots ,m\}$. Then the solution of the recurrence relation 
			\eqref{eq:MinusRecur} and \eqref{eq:PlusRecur} is given by 
      \begin{align}
         \omega^{(-)}_{2m}=&\mathcal{T}^{(-)}_m(\lambda)\omega^{(+)}_0,\nonumber\\
         \omega^{(+)}_{2m}=& \mathcal{T}^{(+)}_m(\lambda)\omega^{(+)}_0,
      \end{align}
     with the boundary data $\omega^{(-)}_0=0$, $\omega^{(+)}_0=\varphi$.
		
		 If $\lambda = \tfrac{\beta}{2}-N$ for some $N\in\{1,\ldots ,m\}$, this 
		 solution holds for all $m<N$.
    \end{theorem}
    \begin{bew}
      In order to shorten the notation, we introduce 
      \begin{align*}
         B^{(+)}_m&:=\big[(\lambda-\beta)\prod_{k=1}^{m}a_{2k}\big]^{-1}. 
      \end{align*}
      First of all, we verify the recurrence \eqref{eq:PlusRecur}. 
      In terms of solution operators, it reads
      \begin{align*}
        a_{2m}\mathcal{T}^{(+)}_m(\lambda)&=b_{2m}(\tfrac{J}{2n})\mathcal{T}^{(+)}_{m-1}(\lambda)
         +c_{2m}(\tfrac{J}{2n})^2\mathcal{T}^{(+)}_{m-2}(\lambda)
          +\Delta^h\mathcal{T}^{(+)}_{m-1}(\lambda)\\
        &\quad+2\d\mathcal{T}^{(-)}_m(\lambda) -2(\tfrac{J}{2n})^2\d \mathcal{T}^{(-)}_{m-2}(\lambda).
      \end{align*}
      This can be decomposed, due to $\d^2=0=(\delta^h)^2$ and 
      dealing with $\d\delta^h$ and $\delta^h \d$ as independent 
			commuting variables, into two independent claims. The first is
      \begin{align}
         a_{2m}B^{(+)}_ms^{(+)}_m(y^{(+)})
           =B^{(+)}_{m-1}\big[\tfrac{n}{2J}\delta^h \d+\tfrac 14 b_{2m}\big]s^{(+)}_{m-1}(y^{(+)})
             +\tfrac{1}{16} c_{2m}B^{(+)}_{m-2}s^{(+)}_{m-2}(y^{(+)}),\label{eq:ScalarVersion}
      \end{align}
      while the second is given by 
      \begin{align}\label{eq:FormVersion}
         a_{2m}B^{(+)}_m s^{(1)}_m(y^{(1)})&=
           B^{(+)}_{m-1}\big[\tfrac{n}{2J}\d\delta^h+\tfrac 14 b_{2m} \big]s^{(1)}_{m-1}(y^{(1)})
           +\tfrac{1}{16}c_{2m}B^{(+)}_{m-2} s^{(1)}_{m-2}(y^{(1)})\notag\\
         &\quad+B^{(+)}_{m-1}(\lambda-\beta) \tfrac{n}{2J} \d\delta^h  s^{(+)}_{m-1}(y^{(+)}) \notag\\
         &\quad+2B^{(+)}_{m-1}\tfrac{n}{2J}\d s^{(-)}_m(y^{(-)}) \delta^h
           -\tfrac 18 B^{(+)}_{m-3} \tfrac{n}{2J} \d s^{(-)}_{m-2}(y^{(-)}) \delta^h .
      \end{align}
      Note that 
      \begin{align*}
        a_{2k}B^{(+)}_k&=B^{(+)}_{k-1},\\
        \tfrac 14 b_{2m}&=2(m-1)(\lambda+m-1)+\tfrac{\beta}{2}\lambda,\\
        \tfrac{1}{16}c_{2m}a_{2m-2}&=-(m-1)(\lambda+m-2)(\lambda-\tfrac{\beta}{2}+m-1)(\tfrac{\beta}{2}+m-2). 
      \end{align*}
      We first notice that Equation \eqref{eq:ScalarVersion} was proved 
			in Proposition \ref{Identities4}. 
      Now we proceed with Equation \eqref{eq:FormVersion}. By $(\delta^h)^2=0$ and 
			Lemma A.1, we have
      \begin{align*}
         \d\delta^h s^{(+)}_{k}(y^{(+)})
            &=s^{(+)}_k\big(\tfrac{\beta}{2}(\tfrac{\beta}{2}-1)\big)  \d\delta^h
            =(\lambda)_{k}(\tfrac{\beta}{2})_{k} \d\delta^h .
      \end{align*}
      Furthermore,  
      \begin{align*}
        \d s^{(-)}_k(y^{(-)})\delta^h =\d\delta^h s^{(-)}_k(y^{(1)})
      \end{align*}
      for all $k\in\N$. Applying Theorem \ref{Identities6} 
      to $s^{(1)}_k(y^{(1)})$, for $k=m, m-1,m-2$, allows to rewrite Equation 
      \eqref{eq:FormVersion} just in terms of $s^{(-)}_k(y^{(1)})$ for appropriate
			collection of values of $k$. It turns out that 
      Equation \eqref{eq:FormVersion} is equivalent to the three-times repeated 
			application of 
			the recurrence relation in Proposition \ref{Identities3} to 
      \begin{multline*}
         (\lambda-\beta+2m)s^{(-)}_{m+1}(y^{(1)})
          -2(m-1)(\lambda+2m-2)(\lambda-\tfrac{\beta}{2}+m-1)s^{(-)}_m(y^{(1)})\\
          +(m-1)(m-2)(\lambda-\beta+2m-4)(\lambda-\tfrac{\beta}{2}+m-2)_2s^{(-)}_{m-1}(y^{(1)}),
      \end{multline*}
      which finally proves Equation \eqref{eq:FormVersion}. 

      Now we proceed to prove  
      \begin{align}
         A_{2m}\mathcal{T}^{(-)}_m&(\lambda)=B_{2m}(\tfrac{J}{2n})\mathcal{T}^{(-)}_{m-1}(\lambda)
            +C_{2m}(\tfrac{J}{2n})^2 \mathcal{T}^{(-)}_{m-2}(\lambda)
            +D_{2m}(\tfrac{J}{2n})^3 \mathcal{T}^{(-)}_{m-3}(\lambda)\notag\\ 
        &\quad+\Delta^h \mathcal{T}^{(-)}_{m-1}(\lambda)-(\tfrac{J}{2n})\Delta^h \mathcal{T}^{(-)}_{m-2}(\lambda)
            +2\delta^h \mathcal{T}^{(+)}_{m-1}(\lambda)
            +2(\tfrac{J}{2n}) \delta^h \mathcal{T}^{(+)}_{m-2}(\lambda). \label{eq:NegFormVersion}
      \end{align}
      Using two ingredients: due to $(\delta^h)^2=0$, we have
      \begin{align*}
         \Delta^h \mathcal{T}^{(-)}_k(\lambda)&=\delta^h d \mathcal{T}^{(-)}_k(\lambda),\\
         \delta^h s^{(1)}_k(y^{(1)})&=s^{(1)}_k(y^{(-)})\delta^h,
      \end{align*}
      and due to Lemma A.1., we get
      \begin{align*}
        \delta^h s^{(+)}_k(y^{(+)})=\sum_{j=0}^k C^{(+)}_j(k)(\tfrac{\beta}{2}-j)_{2j}\delta^h
          = (\lambda)_k(\tfrac{\beta}{2})_k\delta^h,
      \end{align*}
      we see that Equation \eqref{eq:NegFormVersion} is equivalent to 
      \begin{align}\label{eq:FinalRecu}
         A_{2m}B^{(+)}_{m-1} s^{(-)}_m(y^{(-)})&=B^{(+)}_{m-2}\big[\tfrac{n}{2J}\delta^h d
            +\tfrac 14 B_{2m}\big] s^{(-)}_{m-1}(y^{(-)})\notag\\
         &\quad-\tfrac 14 B^{(+)}_{m-3}\big[\tfrac{n}{2J}\delta^h d
            -\tfrac 14 C_{2m}\big] s^{(-)}_{m-2}(y^{(-)})\notag\\
        &\quad+\tfrac{1}{64}B^{(+)}_{m-4} D_{2m}s^{(-)}_{m-3}(y^{(-)})\notag\\
        &\quad+2B^{(+)}_{m-1}s^{(1)}_{m-1}(y^{(-)})
          +2B^{(+)}_{m-1}(\lambda)_{m-1}(\lambda-\beta)(\tfrac{\beta}{2})_{m-1}\notag\\
        &\quad+\tfrac 12 B^{(+)}_{m-2}s^{(1)}_{m-2}(y^{(-)})
           +\tfrac 12 B^{(+)}_{m-2}(\lambda)_{m-2}(\lambda-\beta)(\tfrac{\beta}{2})_{m-2}.
      \end{align}
      We replace the terms $s^{(1)}_{m-1}(y^{(-)})$ and $s^{(1)}_{m-2}(y^{(-)})$, 
      using Theorem \ref{Identities6}, by $s^{(-)}_k(y^{(-)})$ for appropriate 
			collection of values of $k$. 
      Then it turns out that Equation \eqref{eq:FinalRecu} is equivalent to the 
			two-times application of the recurrence relation in 
      Proposition \ref{Identities3} to 
      \begin{align*}
         s^{(-)}_m(y^{(-)})-a_{2m-4}s^{(-)}_{m-1}(y^{(-)}).
      \end{align*}
      This completes the proof. 
    \end{bew}


\section{Applications: Branson-Gover operators on Einstein manifolds}\label{BransonGoverOperators}
  This section is focused on the origin and properties of the 
	Branson-Gover operators and their derived quantities on Einstein manifolds. 
  In addition, we present another proof of a result in \cite{GoverSilhan}
	on the decomposition of Branson-Gover operators 
  as a product of second-order differential operators.

  Let $(M,h)$ be a Riemannian manifold of dimension $n$. 
  For $p=0,\ldots, n$ the Branson-Gover operators \cite{BransonGover} are differential 
  operators 
  \begin{align*}
    L_{2N}^{(p)}:\Omega^p(M)\to\Omega^p(M)
  \end{align*}
  of order $2N$, for $N\in\N$ ($N\leq \tfrac n2$ for even $n$), of the form 
  \begin{align*}
    L_{2N}^{(p)}=(\tfrac{n-2p}{2}+N)(\delta^h\d)^N+(\tfrac{n-2p}{2}-N)(\d\delta^h)^N+
		\mathrm{LOT},
  \end{align*} 
  where $\mathrm{LOT}$ is the shorthand notation for the 
  lower order (curvature correction) terms. They generalize 
  the GJMS operator \cite{GJMS}
  \begin{align*}
    P_{2N}=(\Delta^h)^N+\mathrm{LOT}:\mathcal{C}^\infty(M)\to\mathcal{C}^\infty(M)
  \end{align*}
  in the sense that $L_{2N}^{(0)}=(\tfrac n2+N)P_{2N}$. The key property of 
	Branson-Gover operators is that they are conformally covariant,  
  \begin{align*}
    e^{(\tfrac n2 -p+N)\sigma} \circ\hat{L}_{2N}^{(p)}= L_{2N}^{(p)}\circ e^{(\tfrac n2 -p-N)\sigma},\quad \sigma\in \mathcal{C}^\infty(M) .
  \end{align*}
  Here $\hat{\cdot}$ denotes the evaluation with 
	respect to the conformally related metric 
  $\hat{h}=e^{2\sigma}h$. In the case of even dimensions $n$ 
	and $p\leq \tfrac n2$, the critical Branson-Gover operators factorize 
  \begin{align} 
    L_{n-2p}^{(p)}&=G_{n-2p-1}^{(p+1)}\circ \nonumber\d\\
    &=\delta^h\circ Q_{n-2p-2}^{(p+1)}\circ \d    \label{L-Q-relation}
  \end{align}
  by two additional differential operators 
  \begin{align*}
   G_{n-2p-1}^{(p+1)} = (n-2p)\delta^h\circ\big( (\d\delta)^{\tfrac n2-p-1}+\mathrm{LOT}\big)
     & :\Omega^{p+1}(M)\to\Omega^p(M),\\
   Q_{n-2p-2}^{(p+1)}= (n-2p)\big((\d\delta^h)^{\tfrac n2-p-1}+\mathrm{LOT}\big) &:\Omega^{p+1}(M)\to\Omega^{p+1}(M),
  \end{align*}
  called the gauge companion and the $Q$-curvature operator, 
  respectively. Similarly to $L_{2N}^{(p)}$, these relatives are quite complicated operators
  in general, but in the case when the underlying metric is flat 
  or Einstein we shall present closed formulas for them.  

  Now let $(M,h)$ be an $n$-dimensional Einstein manifold with normalization given by 
  $Ric(h)=\tfrac{2J}{n}(n-1)h$ for a constant $J\in\R$. 
  By \cite{AG}, it follows that one can recover Branson-Gover 
  operators as residues of solution operators, 
  see Equation \eqref{eq:SolutionOperators}. 
  More precisely, we have 
  \begin{align}\label{eq:BransonGoverOp}
      \Res_{\lambda=\tfrac{\beta}{2}-N}\big(\mathcal{T}^{(+)}_{2N}(\lambda)\big)
        &\sim (\tfrac{2J}{n})^N\Big[(\tfrac{\beta}{2}+N)R_N(y^{(+)};0)+(\tfrac{\beta}{2}-N)R^{(1)}_N(y^{(1)})\Big].
  \end{align}
  The right hand side on the previous display is just the Branson-Gover operator of order 
  $2N$, $N\in\N$, acting on differential $p$-forms:
  \begin{align}\label{eq:BGDefinition}
    L_{2N}^{(p)}=(\tfrac{2J}{n})^N\Big[(\tfrac{\beta}{2}+N)R_N(y^{(+)};0)+(\tfrac{\beta}{2}-N)R^{(1)}_N(y^{(1)})\Big].
  \end{align}
  Note that there is no obstruction for its existence 
  in even dimensions $n$. The normalization 
  \begin{align}
    \bar{L}_{2N}^{(p)}:&=(\tfrac{\beta}{2}-N+1)_{2N-1} L_{2N}^{(p)}\notag\\
    &=(\tfrac{2J}{n})^N\Big[(\tfrac{\beta}{2}-N+1)_{2N }R_N(y^{(+)};0)+(\tfrac{\beta}{2}-N)_{2N}R^{(1)}_N(y^{(1)})\Big]
  \end{align}
  has the effect that the factors appearing in Theorem \ref{MainTheorem} 
  are differential operators with polynomial coefficients. 
  \begin{bem}\label{ExceptionalValues}
    The normalization factor $(\tfrac{\beta}{2}-N+1)_{2N-1}$ can vanish only 
		in even dimensions $n$, due to $\beta=n-2p$. 
    The zero locus is characterized by $l\in \{0,\ldots,N-1\}$ such that 
    $\tfrac{\beta}{2}=l$, or by the existence of $l\in\{1,\ldots N-1\}$ 
		such that $\tfrac{\beta}{2}=-l$. 
    The last case can be excluded by choosing $p\leq \tfrac n2$.   
  \end{bem}
  \begin{sat}\label{BransonGoverRecurrence}
    Let $N\in\N$ and $p=0,\ldots,n$ when $n$ is odd, and $N\in\N$ and $p=0,\ldots,n$ 
    such that $(\tfrac{\beta}{2}-N+1)_{2N-1}\neq 0$ 
    when $n$ is even. 
    The normalized Branson-Gover operators $\bar{L}_{2N}^{(p)}$
    satisfy the recurrence relation 
    \begin{align}
         \bar{L}_{2N}^{(p)}=\big[& (\tfrac{\beta}{2}+N)(\tfrac{\beta}{2}-N+1) \delta^h\d
         +(\tfrac{\beta}{2}+N-1)(\tfrac{\beta}{2}-N)\d\delta^h\notag\\
     &+(\tfrac{\beta}{2}-N)(\tfrac{\beta}{2}-N+1)(\tfrac{\beta}{2}+N-1)(\tfrac{\beta}{2}+N)\tfrac{2J}{n}\big]\bar{L}_{2N-2}^{(p)},\label{eq:BGRec}
    \end{align}
    for $\bar{L}^{(p)}_0:=\id$. 
  \end{sat}
  \begin{bew}
    We shall use 
    \begin{align*}
      \tfrac{2J}{n}R_m(y^{(+)};0)&=\big[\delta^h\d+(\tfrac{\beta}{2}+m-1)(\tfrac{\beta}{2}-m)\tfrac{2J}{n}\big] 
        R_{m-1}(y^{(+)};0),\\
      \tfrac{2J}{n}R^{(1)}_m(y^{(1)})&=\big[\d\delta^h+(\tfrac{\beta}{2}+m)(\tfrac{\beta}{2}-m+1)\tfrac{2J}{n}\big]
        R^{(1)}_{m-1}(y^{(1)})+ (\tfrac{\beta}{2}-m+2)_{2m-2} \d\delta^h,
    \end{align*}
    where the first equality is easily verified and 
		the second follows from Theorem \ref{Identities2}.
    In addition, we need an elementary identity 
    \begin{align*}
      \d\delta^h R_{N-1}(y^{(+)};0)=(\tfrac{\beta}{2}-N+1)_{2N-2}\d\delta^h.
    \end{align*}
    We start with the evaluation of the right hand side of Equation \eqref{eq:BGRec}. We have 
    \begin{align*}
      \Big[& (\tfrac{\beta}{2}+N)(\tfrac{\beta}{2}-N+1) \delta^h\d
         +(\tfrac{\beta}{2}+N-1)(\tfrac{\beta}{2}-N)\d\delta^h \\
     &\quad+(\tfrac{\beta}{2}-N)(\tfrac{\beta}{2}-N+1)(\tfrac{\beta}{2}+N-1)(\tfrac{\beta}{2}+N)\tfrac{2J}{n}\Big]\times\\
     &\quad \times  (\tfrac{2J}{n})^{N-1}\Big[(\tfrac{\beta}{2}-N+2)_{2N-2 }R_{N-1}(y^{(+)};0)
        +(\tfrac{\beta}{2}-N+1)_{2N-2}R^{(1)}_{N-1}(y^{(1)})\Big]\\
     &=(\tfrac{2J}{n})^{N-1}(\tfrac{\beta}{2}-N+1)_{2N}\big[\delta^h\d
       +(\tfrac{\beta}{2}-N)(\tfrac{\beta}{2}+N-1)\tfrac{2J}{n}\big]R_{N-1}(y^{(+)};0)\\
     &\quad+(\tfrac{2J}{n})^{N-1}(\tfrac{\beta}{2}-N)_{2N}\big[\d\delta^h
       +(\tfrac{\beta}{2}-N+1)(\tfrac{\beta}{2}+N)\tfrac{2J}{n}\big]R^{(1)}_{N-1}(y^{(1)})\\
     &\quad+(\tfrac{2J}{n})^{N-1}(\tfrac{\beta}{2}-N)(\tfrac{\beta}{2}+N-1)(\tfrac{\beta}{2}-N+2)_{2N-2}\d\delta^h R_{N-1}(y^{(+)};0).
    \end{align*}
    The preparatory identities above ensure that this equals to 
    \begin{align*}
      (\tfrac{2J}{n})^N\Big[(\tfrac{\beta}{2}-N+1)_{2N}R_{N}(y^{(+)};0)+(\tfrac{\beta}{2}-N)_{2N}R^{(1)}_N(y^{(1)})\Big]
    \end{align*}
    and the proof is complete. 
  \end{bew}
 
  The recurrence relation for $\bar{L}_{2N}^{(p)}$ implies part of the 
	result \cite[Theorem $5.3$]{GoverSilhan}.
  \begin{theorem}\label{MainTheorem}
    Let $N\in\N$ and $p=0,\ldots,n$ when $n$ is odd, and $N\in\N$ and $p=0,\ldots,n$ 
    such that $(\tfrac{\beta}{2}-N+1)_{2N-1}\neq 0$ 
    when $n$ is even. The normalized Branson-Gover operators 
    $\bar{L}_{2N}^{(p)}$ factorize as
    \begin{align}
      \bar{L}_{2N}^{(p)}=\prod_{l=1}^N\big[ &(\tfrac{\beta}{2}+N-l+1)(\tfrac{\beta}{2}-N+l)\delta^h\d
         +(\tfrac{\beta}{2}+N-l)(\tfrac{\beta}{2}-N+l-1)\d\delta^h \nonumber \\
      &+(\tfrac{\beta}{2}-N+l-1)(\tfrac{\beta}{2}-N+l)(\tfrac{\beta}{2}+N-l)(\tfrac{\beta}{2}+N-l+1)\tfrac{2J}{n}\big].
    \end{align}
  \end{theorem}
  In the setting of Theorem \ref{MainTheorem} we have 
  \begin{align}
    L_{2N}^{(p)}=\frac{\beta}{2}\prod_{k=1}^N \big[ \tfrac{\beta+2k}{\beta+2k-2}\delta^h\d
      +\tfrac{\beta-2k}{\beta-2k+2}\d\delta^h+(\tfrac{\beta}{2}-k)(\tfrac{\beta}{2}+k)\tfrac{2J}{n} 
			\big],
  \end{align}
  which holds in even dimensions $n$ only for $\beta\notin\{0,2,\ldots,2N-2\}$.
  
  Now we discuss the exceptional cases when 
	$(\tfrac{\beta}{2}-N+1)_{2N-1}=0$ and thus completing the result \cite[Theorem $5.3$]{GoverSilhan}. By Remark \ref{ExceptionalValues} 
  we have to consider $\beta=2l$ for $l\in\{0,\ldots,N-1\}$. For these values 
  the polynomials $R_N(y^{(+)};0)$ and $R^{(1)}_N(y^{(1)})$ factorize by 
	$\delta^h\d$ and $\d\delta^h$, respectively, and this influences 
	the factorization of $L_{2N}^{(p)}$. 

  \begin{theorem}\label{BGForExceptionals}
		Let $n$ be even, $N\in\N$ and $p\leq \tfrac n2$. 

    \begin{enumerate}
    \item Let $\beta=2l$ for some $l\in\{1,\ldots,N-1\}$. The Branson-Gover operator 
		of order $2N$ factorizes by 
     \begin{align}
        L_{2N}^{(p)}=-\tfrac{2l}{2l-1}\widetilde{P}^{(p)}\times 
         \prod_{k=1,k\neq l,l+1}^{N}\big[\tfrac{l+k}{l+k-1}\delta^h\d
         +\tfrac{l-k}{l-k+1}\d\delta^h+(l-k)(l+k)\tfrac{2J}{n}\big]  ,
     \end{align}
     where 
     \begin{align*}
        \tilde{P}^{(p)}:=\big[\tfrac{2l+1}{2}\delta^h \d+\tfrac{2l-1}{2}\d\delta^h\big]
      \big[\d\delta^h-\delta^h\d+2l\tfrac{2J}{n}\big]
     \end{align*}
     is a fourth-order differential operator. 
     \item Let $\beta=0$. The Branson-Gover operator of order $2N$ factorizes by
     \begin{align}
       L_{2N}^{(p)}=N\big[\delta^h\d-\d\delta^h\big]\times
         \prod_{k=2}^{N}\big[\delta^h\d
         +\d\delta^h-k(k-1)\tfrac{2J}{n}\big].
     \end{align}
     \end{enumerate}
  \end{theorem}
  \begin{bew}
    First note that $\tilde{P}^{(p)}$ decomposes, due to $\d^2=0=(\delta^h)^2$, as
    \begin{align*}
      \tilde{P}^{(p)}=-\tfrac{2l+1}{2}\delta^h \d\big[\delta^h\d-2l \tfrac{2J}{n}\big]  
       +\tfrac{2l-1}{2}\d\delta^h\big[\d\delta^h+2l\tfrac{2J}{n}\big]
      =:\widetilde{P}^{(p)}_1+\widetilde{P}^{(p)}_2.
    \end{align*}
    For $l\in\{1,\ldots,N-1\}$ such that $\beta=2l$ the polynomials $R_m(y^{(+)};0)$ 
		and $R^{(1)}_m(y^{(1)})$, cf. Equation \eqref{eq:Polynomal1} and \eqref{eq:Polynomal2}, 
		satisfy by Theorem \ref{Identities1}
    \begin{align*}
       R_m(y^{(+)};0)&=-\tfrac{2}{2l+1} (\tfrac{n}{2J})^2\widetilde{P}^{(p)}_1\times
         \prod_{k=1,k\neq l,l+1}^{N}\big[\tfrac{n}{2J}\delta^h \d+(l-k)(l+k-1) \big] ,\\
      R^{(1)}_m(y^{(1)})&=R_m(y^{(1)};0)+0\\
      &=\tfrac{2}{2l-1}(\tfrac{n}{2J})^2\widetilde{P}^{(p)}_2\prod_{k=1,k\neq l,l+1}^{N}\big[\tfrac{n}{2J}\d \delta^h
        +(l+k)(l-k+1) \big].
    \end{align*}
    Since $\widetilde{P}^{(p)}_1$ and $\widetilde{P}^{(p)}_2$ factorize by $\delta^h\d$ and $\d\delta^h$, 
    respectively, the last display is equivalent to 
    \begin{align*}
       R_m(y^{(+)};0)&=-\tfrac{2l}{2l-1}\tfrac{1}{l+N}(\tfrac{n}{2J})^N\widetilde{P}^{(p)}_1\times\notag\\
       &\quad\quad\times \prod_{k=1,k\neq l,l+1}^{N}\big[\tfrac{l+k}{l+k-1}\delta^h \d+\tfrac{l-k}{l-k+1}\d\delta^h
         +(l-k)(l+k)\tfrac{2J}{n} \big],\\
      R^{(1)}_m(y^{(1)})&=-\tfrac{2l}{2l-1}\tfrac{1}{l-N}(\tfrac{n}{2J})^N \widetilde{P}^{(p)}_2\times\notag\\
      &\quad\quad\times  \prod_{k=1,k\neq l,l+1}^{N}\big[\tfrac{l-k}{l-k+1}\d \delta^h+\tfrac{l+k}{l+k-1}\delta^h\d
        +(l-k)(l+k)\tfrac{2J}{n} \big].
    \end{align*}
    Hence, the definition of $L_{2N}^{(p)}$, see \eqref{eq:BGDefinition}, gives 
    \begin{align*}
      L_{2N}^{(p)}=-\tfrac{2l}{2l-1}\widetilde{P}^{(p)}\times 
       \prod_{k=1,k\neq l,l+1}^{N}\big[\tfrac{l-k}{l-k+1}\d \delta^h+\tfrac{l+k}{l+k-1}\delta^h\d
        +(l-k)(l+k)\tfrac{2J}{n} \big],
    \end{align*}
    which proves the first claim. 

    For $\beta=0$, the polynomials $R_m(y^{(+)};0)$ and $R^{(1)}_m(y^{(1)})$ 
    (cf. Equation \eqref{eq:Polynomal1} and \eqref{eq:Polynomal2}) satisfy 
		by Theorem \ref{Identities1}
    \begin{align*}
       R_m(y^{(+)};0)&=\tfrac{n}{2J}\delta^h\d\prod_{k=2}^{N}\big[\tfrac{n}{2J}\delta^h \d
         -k(k-1) \big],\\
      R^{(1)}_m(y^{(1)})&=R_m(y^{(1)};0)+0\notag\\
      &=\tfrac{n}{2J}\d \delta^h\prod_{k=2}^{N}\big[\tfrac{n}{2J}\d \delta^h
        -k(k-1) \big].
    \end{align*}
    Since they factorize by $\delta^h\d$ and $\d\delta^h$, respectively, we can write 
    \begin{align*}
      R_m(y^{(+)};0)&=(\tfrac{n}{2J})^N\delta^h\d\prod_{k=2}^{N}\big[\delta^h \d+\d\delta^h
         -k(k-1)\tfrac{2J}{n} \big],\\
      R^{(1)}_m(y^{(1)})&=(\tfrac{n}{2J})^N\d \delta^h\prod_{k=2}^{N}\big[\d \delta^h+\delta^h\d
        -k(k-1)\tfrac{2J}{n} \big].
    \end{align*}
    Hence the definition of $L_{2N}^{(p)}$, see \eqref{eq:BGDefinition}, gives 
    \begin{align*}
      L_{2N}^{(p)}&=N\big[\delta^h\d-\d\delta^h\big]\times
        \prod_{k=2}^N\big[\delta^h\d+\d\delta^h-k(k-1)\tfrac{2J}{n}  \big]
    \end{align*}
    and the proof is complete.
  \end{bew}
  
  \begin{bem}
    Although the results in Theorems \ref{MainTheorem} and \ref{BGForExceptionals} are not new we think that we could enlight 
    the appearance of the factorizing structure of the Branson-Gover operators from a purely combinatorical point of view. 
  \end{bem}

  
  From now on let $n$ be even. We proceed with explicit formulas 
	for the critical Branson-Gover operator, gauge companion operator 
	and $Q$-curvature operator on Einstein manifolds. 
  \begin{theorem}
    The critical Branson-Gover operator $L_{n-2p}^{(p)}$ is given by the product formula
    \begin{align}
      L_{n-2p}^{(p)}= (n-2p) \delta^h\d \circ\prod_{l=1}^{\tfrac{n-2p}{2}-1}
        \big[\delta^h\d+(\tfrac{\beta}{2}-l)(\tfrac{\beta}{2}+l-1)(\tfrac{2J}{n})\big].
    \end{align}
  \end{theorem}
  \begin{bew}
    It follows from Equation \eqref{eq:BransonGoverOp} that 
    \begin{align*}
      L_{n-2p}^{(p)}&=(n-2p)(\tfrac{2J}{n})^{\tfrac{n-2p}{2}}  R_{\tfrac{n-2p}{2}}(y^{(+)};0)\\
      &=(n-2p) \prod_{l=1}^{\tfrac{n-2p}{2}}\big[\delta^h\d+(\tfrac{\beta}{2}-l)(\tfrac{\beta}{2}+l-1)\tfrac{2J}{n}\big].
    \end{align*}
    Note that the last factor in that product reduces to 
    \begin{align*}
      \delta^h\d+(\tfrac{\beta}{2}-\tfrac{n-2p}{2})(\tfrac{\beta}{2}+\tfrac{n-2p}{2}-1)\tfrac{2J}{n}=\delta^h\d,
    \end{align*}
    since $\beta=n-2p$. This completes the proof. 
  \end{bew}

  Consequently by \eqref{L-Q-relation}, we found the explicit formulas for the $Q$-curvature operator 
  \begin{align*}
    Q^{(p+1)}_{n-2p-2}=(n-2p)\prod_{l=1}^{\tfrac{n-2p}{2}-1}\big[\d\delta^h+(\tfrac{\beta}{2}-l)(\tfrac{\beta}{2}+l-1)(\tfrac{2J}{n})\big]
  \end{align*}
  and the gauge companion operator 
  \begin{align*}
    G^{(p+1)}_{n-2p-1}=(n-2p)\delta^h\circ \prod_{l=1}^{\tfrac{n-2p}{2}-1}\big[\d\delta^h+(\tfrac{\beta}{2}-l)(\tfrac{\beta}{2}+l-1)(\tfrac{2J}{n})\big].
  \end{align*}
  Obviously, 
  \begin{align*}
     L_{n-2p}^{(p)}&=G^{(p+1)}_{n-2p-1}\circ\d=\delta^h\circ  Q^{(p+1)}_{n-2p-2}\circ \d,
  \end{align*}
  which is the famous double factorization of the critical Branson-Gover operator. 

  \begin{bem}
    Let $(M,h)=(\R^n,\langle\cdot,\cdot\rangle)$ be the Euclidean space. The 
		explicit formulas for $L_{2N}^{(p)}$, 
		$Q^{(p+1)}_{n-2p-2}$ and $G^{(p+1)}_{n-2p-1}$ immediately imply after 
		setting $J=0$ that 
    \begin{align*}
      L_{2N}^{(p)}&=(\tfrac{n-2p}{2}+N)(\delta^h\d)^N+(\tfrac{n-2p}{2}-N)(\d\delta^h)^N,\\
      Q^{(p+1)}_{n-2p-2}&=(n-2p)(\d\delta^h)^{\tfrac{n-2p}{2}-1},\\
      G^{(p+1)}_{n-2p-1}&=(n-2p)\delta^h (\d\delta^h)^{\tfrac{n-2p}{2}-1}.
    \end{align*}
  \end{bem}

\appendix
\section{Generalized Hypergeometric functions and dual Hahn polynomials}\label{AppendixA}

Here we summarize a few basic conventions and definitions related to
generalized hypergeometric functions.

The Pochhammer symbol of a complex number $a\in\C$ is defined 
by $(a)_m:=\frac{\Gamma(a+m)}{\Gamma(a)}$, for $m\in\Z$ and $(a)_0:=1$. If $m$ is a positive 
integer, $(a)_m=a(a+1)\cdots(a+m-1)$. 

The generalized hypergeometric function ${}_pF_q$ of type $(p,q)$,
$p,q\in\N$, is defined by
\begin{align}
    {}_pF_q\left[\begin{matrix}a_1\;,\;\ldots\;,\;a_p\\b_1\;,\;\ldots\;,\;b_q\end{matrix};z\right]
      :=\sum_{l=0}^\infty\frac{(a_1)_l\ldots (a_p)_l}{(b_1)_l\ldots (b_q)_l}\frac{z^l}{l!},
\end{align}
for $a_i\in\C$ ($1\leq i\leq q$), $b_j\in\C\setminus\{-\N_0\}$ ($1\leq j\leq q$), and $z\in\C$.

The dual Hahn polynomials \cite{KMcG} are defined by 
\begin{align*}
  R_m\big(\lambda(n);a,b,N\big)={}_3F_2\left[\begin{matrix}-m\;,\;-n\;,\;n+a+b+1\\a+1\;,\;-N+1\end{matrix};1\right]
\end{align*}
for $n,N\in\N$, $m=0,\ldots,N-1$, $a,b\in\C$ such that the real parts fulfill 
$\mathfrak{R}(a),\mathfrak{R}(b)>-1$ and $\lambda(n):=n(n+a+b+1)$.
 
Here we prove the two identities for collections of
$C^{(+)}_k(m)$ and $ D^{(1)}_k(m)$, defined in 
Equations \eqref{eq:Def-s+} and \eqref{eq:Def-s1}:
\begin{lem1}
      We have 
      \begin{align*}
        \sum_{k=0}^{m} C^{(+)}_k(m) (\tfrac{\beta}{2}-k)_{2k}&=(\lambda)_m (\tfrac{\beta}{2})_m,\\
        \sum_{k=0}^{m} D^{(1)}_k(m) (\tfrac{\beta}{2}-k+1)_{2k}&=(\lambda)_m (\tfrac{\beta}{2})_{m}(\lambda-\beta).
      \end{align*}
    \end{lem1}
    \begin{bew}
      By standard Pochhammer identities, we rewrite the sums as 
      \begin{align*}
        \sum_{k=0}^{m}& C^{(+)}_k(m) (\tfrac{\beta}{2}-k)_{2k}
          =(\lambda-\tfrac{\beta}{2}+1)_m(\tfrac{\beta}{2})_m
          \sum_{k=0}^m\frac{(-m)_k (-\tfrac{\beta}{2}+1)_k}{k! (\lambda-\tfrac{\beta}{2}+1)_k}\\
        &=(\lambda-\tfrac{\beta}{2}+1)_m(\tfrac{\beta}{2})_m 
           \times {}_2F_1(-m, -\tfrac{\beta}{2}+1; \lambda-\tfrac{\beta}{2}+1;1),\\
        \sum_{k=0}^{m}& D^{(1)}_k(m) (\tfrac{\beta}{2}-k+1)_{2k}
          =(\lambda-\tfrac{\beta}{2}+1)_m(\tfrac{\beta}{2}+1)_{m-1}
          \sum_{k=0}^m\frac{(-m)_k (-\tfrac{\beta}{2})_k}{k! (\lambda-\tfrac{\beta}{2}+1)_k}\times\\
        &\quad\quad\times \big[ k(\lambda+\beta+2m)+\tfrac{\beta}{2}(\lambda-\beta-2m) \big]\\
        &=(\lambda-\tfrac{\beta}{2}+1)_m(\tfrac{\beta}{2}+1)_{m-1}
          \Big[\tfrac{m \tfrac{\beta}{2}(\lambda+\beta+2m)}{\lambda-\tfrac{\beta}{2}+1}
          \times {}_2F_1(-m+1, -\tfrac{\beta}{2}+1; \lambda-\tfrac{\beta}{2}+2;1)\\
        &\quad\quad+\tfrac{\beta}{2}(\lambda-\beta-2m)
            \times {}_2F_1(-m,-\tfrac{\beta}{2}; \lambda-\tfrac{\beta}{2}+1;1)\Big]. 
      \end{align*}
      Applying the Chu-Vandermonde identity, 
      \begin{align*}
        {}_2F_1\left[\begin{matrix}-m,\alpha\\ \gamma  \end{matrix};1\right]=\frac{(\gamma-\alpha)_m}{(\gamma)_m}
      \end{align*} 
      with $m\in \N$ and appropriate $\alpha,\gamma\in\C$, the first claim follows. 
      The proof of the second statement is based on the identity 
      \begin{align*}
        m(\lambda+\beta+2m)(\lambda+1)_{m-1}+(\lambda-\beta-2m)(\lambda+1)_m=(\lambda)_m(\lambda-\beta).
      \end{align*} 
      This completes the proof. 
    \end{bew}

\section{Poincar\'e-Einstein metric construction}\label{AppendixB}
Here we briefly review the content of Poincar\'e-Einstein metric construction, \cite{FG3}. Let $(M^n,h)$ be an $n$-dimensional 
semi-Riemannian manifold, $n\geq 3$. On $X:=M\times (0,\varepsilon)$ for $\varepsilon>0$, we 
consider the metric
\begin{align*}
  g_+=r^{-2}(\d r^2+h_r),
\end{align*}
for a $1$-parameter family of metrics $h_r$ on $M$ such that $h_0=h$. The requirement of Einstein 
condition on $g_+$ for $n$ odd, 
\begin{align*}
  Ric(g_+)+ng_+=O(r^\infty),
\end{align*}
uniquely determines the family $h_r$, while for $n$ even the conditions 
\begin{align*}
  Ric(g_+)+ng_+=O(r^{n-2}),\\
  \tr(Ric(g_+)+ng_+)=O(r^{n-1}),
\end{align*}
uniquely determine the coefficients $h_{(2)},\ldots,h_{(n-2)}$, $\tilde{h}_{(n)}$ and the trace of $h_{(n)}$ in the 
formal power series
\begin{align*}
  h_r=h+r^2h_{(2)}+\cdots+r^{n-2}h_{(n-2)}+r^n(h_{(n)}+\tilde{h}_{(n)}\log r )+\cdots\quad . 
\end{align*}
For example, we have 
\begin{align*}
  h_{(2)}=-P,\quad h_{(4)}=\frac 14(P^2-\frac{B}{n-4}),
\end{align*}
where $P$ is the Schouten tensor and $B$ is the Bach tensor associated to $h$. 
The metric $g_+$ on $X$ is called Poincar\'e-Einstein metric associated to the 
semi-Riemannian manifold $(M,h)$. 

Two different representatives $h,\widehat{h}\in [h]$ 
in the conformal class lead to Poincar\'e-Einstein 
metrics $g^1_+$ and $g^2_+$ related by a diffeomorphism 
$\Phi:U_1\subset X\to U_2\subset X$, 
where both $U_i$, $i=1,2$, contain $M\times\{0\}$, 
$\Phi|_{M}=\id_M$ and $g_+^1=\Phi^*g^2_+$ (up to a finite 
order in $r$ for even $n$). 

The explicit knowledge of $h_r$ is available in a few cases, 
e.g., \cite{FG3, LeistnerNurowski}. If $(M,h)$ is 
the Euclidean space, one can realize the Poincar\'e-Einstein 
metric as the hyperbolic metric on the upper half space, while 
if $(M,h)$ is an Einstein manifold normalized by 
$Ric(h)=\tfrac{2J}{n}(n-1)h$ for some constant $J\in\R$, 
the $1$-parameter family of metrics $h_r$ is  
\begin{align*}
   h_r=(1-\tfrac{J}{2n}r^2)^2h.
\end{align*}
In both cases, there are no obstructions in even dimensions.

\bibliographystyle{alpha}
\bibliography{BGOnEinstein01052019}

\newcommand{\etalchar}[1]{$^{#1}$}
\begin{thebibliography}{KKM{\etalchar{+}}78}

\bibitem[AG08]{AG}
E.~Aubry and C.~Guillarmou.
\newblock Conformal harmonic forms, {B}ranson-{G}over operators and {D}irichlet
  problem at infinity.
\newblock {\em ArXiv e-prints}, 2008.
\newblock \url{http://arxiv.org/pdf/0808.0552v1.pdf}.

\bibitem[BG05]{BransonGover}
T.~Branson and A.R. Gover.
\newblock Conformally invariant operators, differential forms, cohomology and a
  generalisation of {$Q$-}curvature.
\newblock {\em Communications in Partial Differential Equations},
  30(11):1611--1669, 2005.

\bibitem[Bra93]{Branson5}
T.~Branson.
\newblock The functional determinant.
\newblock {\em Seoul National University Research Institute of Mathematical
  Global Analysis Research Center}, 4:1--101, 1993.

\bibitem[FG11]{FG3}
C.~Fefferman and C.~R. Graham.
\newblock {\em The Ambient Metric ({AM}-178)}.
\newblock Number 178. Princeton University Press, 2011.

\bibitem[Fis13]{Fischmann}
M.~Fischmann.
\newblock {\em Conformally covariant differential operators acting on spinor
  bundles and related conformal covariants}.
\newblock PhD thesis, Humboldt Universit\"at zu Berlin, 2013.
\newblock
  \url{http://edoc.hu-berlin.de/dissertationen/fischmann-matthias-2013-03-04/P%
DF/fischmann.pdf}.

\bibitem[FKS15]{FKS}
M.~Fischmann, C.~Krattenthaler, and P.~Somberg.
\newblock On conformal powers of the {D}irac operator on {E}instein manifolds.
\newblock {\em Mathematische Zeitschrift}, 2015.

\bibitem[GJMS92]{GJMS}
C.R. Graham, R.W. Jenne, L.~Mason, and G.~Sparling.
\newblock Conformally invariant powers of the {L}aplacian, {I}: {E}xistence.
\newblock {\em Journal of the London Mathematical Society}, 2(3):557--565,
  1992.

\bibitem[GLW15]{GoverLatiniWaldron}
R.~Gover, E.~Latini, and A.~Waldron.
\newblock {\em Poincar\'e-{E}instein {H}olography for forms via conformal
  geometry in the bulk}, volume 235.
\newblock American Mathematical Society, 2015.

\bibitem[GMP10]{GMP}
C.~Guillarmou, S.~Moroianu, and J.~Park.
\newblock Eta invariant and {S}elberg zeta function of odd type over convex
  co-compact hyperbolic manifolds.
\newblock {\em Advances in Mathematics}, 225(5):2464--2516, 2010.

\bibitem[G{\v{S}}13]{GoverSilhan}
A.R. Gover and J.~{\v{S}}ilhan.
\newblock Conformally {O}perators on {W}eighted {F}orms; {T}heir
  {D}ecomposition and {N}ull {S}pace on {E}instein {M}anifolds.
\newblock {\em Annales Henri Poincar\'e}, 15(4):679--705, 2013.

\bibitem[GZ03]{GZ}
C.R. Graham and M.~Zworski.
\newblock Scattering matrix in conformal geometry.
\newblock {\em Inventiones mathematicae}, 152:89--118, 2003.

\bibitem[HS01]{HS}
J.~Holland and G.~Sparling.
\newblock Conformally invariant powers of the ambient {D}irac operator.
\newblock {\em ArXiv e-prints}, 2001.
\newblock \url{http://arxiv.org/abs/math/0112033}.

\bibitem[Juh13]{Juhl1}
A.~Juhl.
\newblock Explicit formulas for {GJMS}-operators and {$Q$-}curvatures.
\newblock {\em Geometric and Functional Analysis}, 23:1278--1370, 2013.
\newblock \url{http://arxiv.org/abs/1108.0273}.

\bibitem[KKM{\etalchar{+}}78]{KKMOOT}
M.~Kashiwara, A.~Kowata, K.~Minemura, K.~Okamoto, T.~Oshima, and M.~Tanaka.
\newblock Eigenfunctions of invariant differential operators on a symmetric
  space.
\newblock {\em The Annals of Mathematics}, 107(1):1--39, 1978.

\bibitem[KM61]{KMcG}
S.~Karlin and J.L. McGregor.
\newblock The {H}ahn polynomials, formulas and an application.
\newblock {\em Scripta Mathematica}, 26:33--46, 1961.

\bibitem[LN10]{LeistnerNurowski}
T.~Leistner and P.~Nurowski.
\newblock Ambient {M}etrics for $n$-{D}imensional pp-{W}aves.
\newblock {\em Communications in Mathematical Physics}, 296(3):881--898, 2010.

\end{thebibliography}

\end{document}